\UseRawInputEncoding
\documentclass[10pt]{amsart}
\usepackage{graphicx,amssymb,amsmath}
\newtheorem{thm}{Theorem}[section]
\newtheorem{cor}[thm]{Corollary}
\newtheorem{lem}[thm]{Lemma}
\newtheorem{prop}[thm]{Proposition}

\theoremstyle{definition}

\theoremstyle{remark}
\newtheorem{rem}[thm]{Remark}
\numberwithin{equation}{section}
\usepackage{amsmath,amssymb,amsthm}

\theoremstyle{plain}

\theoremstyle{definition}

\newtheorem*{claim}{Claim}
\newtheorem*{notation}{Notation}

\newtheorem*{main}{Main Theorem}
\newtheorem*{mainc}{Corollary}

\def\Z   {{\bf Z}}

\def\R   {{\bf R}}

\def\P   {{\bf P}}

\def\calK      {{\mathcal K}}

\def\O           {{\mathcal O}}
\def\calE       {{\mathcal E}}

\def\Kahler    {K\"ahler\ }
\def\CY        {Calabi--Yau\ }

\def\Mov       {{\it Mov}}

\def\Movb     {{\overline \Mov}}
\def\Eff         {{\it Eff}}
\def\Effb       {{\overline \Eff}}

\def\vol         {{\rm vol}}

\def\Z   {{\bf Z}}

\def\R   {{\bf R}}

\def\P  {{\bf P}}

\newcommand{\beas}{\begin{eqnarray*}} 
\newcommand{\eeas}{\end{eqnarray*}} 

\begin{document}

\title{Calabi-Yau threefolds with Picard number three}
\author{P.M.H. Wilson} 
\address{Department of Pure Mathematics, University of Cambridge,
16 Wilberforce Road, Cambridge CB3 0WB, UK}
\email {pmhw@dpmms.cam.ac.uk}

\date{10 May 2021}

\begin{abstract}

In this paper, we continue the 
study of boundedness questions for (simply connected) smooth Calabi--Yau threefolds commenced in \cite{WilBd}.  
The diffeomorphism class of such a threefold is known to be determined up to finitely many possibilities by the integral middle cohomology and two integral forms on the integral second cohomology, namely the cubic cup-product form and the linear form given by cup-product with the second chern class.  The question addressed in both papers is whether knowledge of these cubic and linear forms determines the threefold up to finitely many families, that is the moduli of such threefolds is bounded.  If this is true, then in particular the middle integral cohomology would be bounded by knowing these two forms.

Crucial to  this question is the study of rigid non-movable surfaces on the threefold, which are the irreducible surfaces  that deform with any small deformation of the 
complex structure of the threefold but for which no multiple moves on the threefold.   We showed in \cite{WilBd} that 
 if there are no such surfaces, then the answer to the above question is yes.  Moreover if $\rho =2$, the answer was shown to be yes without the further assumption on the Calabi--Yau.
 
 The main results of this paper are for Picard number $\rho =3$, where we prove  boundedness 
 in the case where there is at most one rigid non-movable surface, assuming the cubic form is smooth, thereby defining a real elliptic curve; in the two cases where the Hessian curve is singular, we also assume that 
 the line defined by the second chern class does not intersect this curve at an inflexion point.  The arguments 
 used nicely illustrate the general theory developed by the author in the first half of \cite{WilBd}.
 In addition to the methods described in \cite{WilBd}, a further crucial tool in the proofs will be the classical Steinian involution on the Hessian of an elliptic curve. 
 \end{abstract}


\maketitle

\section*{Introduction}

\rm Throughout this paper, $X$ will denote a (simply connected) smooth complex \CY threefold.  We know that its diffeomorphism class is determined 
up to finitely many possibilities by knowledge of the cup-product cubic form on $H^2 (X, \Z)$ given by $D \mapsto D^3$, the linear form 
on $H^2 (X, \Z)$ given by $D \mapsto D\cdot c_2(X)$ and the middle cohomology $H^3 (X, \Z)$ \cite{Sull, Wall}.  A well-known question is whether the 
\CY threefold $X$ is determined up to finitely many families by the diffeomorphism type.  In the paper \cite{WilBd}, we addressed the 
question  
as to whether $X$ is determined up to finitely many families by the weaker information  of the cubic and linear forms on $H^2 (X, \Z)$, and in particular proved this was true for Picard number $\rho =2$.
In this paper we start the study of higher Picard number.    We saw in \cite{WilBd} that the rigid non-movable surfaces played a central role in this question; if there are no such surfaces 
on $X$, then boundedness was proved in general for all $\rho$.  

Following definitions from \cite{WilBd}, the \it positive index cone \rm of the cubic is the set of classes $L \in H^2 (X, \R) $ for which $L^3 >0$ and the quadratic form given by $D\mapsto L\cdot D^2$ has index $(1, \rho -1)$.  Recall that if $P^\circ$ is the component of the positive index cone on a Calabi--Yau threefold $X$ which contains the \Kahler cone $\calK$, we consider the subcone 
of $P^\circ$ given by the extra conditions that $E\cdot D^2 >0$ for all rigid non-movable surfaces $E$ on $X$, and let $Q$ be the component of this cone which contains the \Kahler cone.  Assuming that there are only finitely many such surfaces $E$ and that $X$ is general in moduli, 
we showed in \cite{WilBd} that any integral class $D$ in $Q$ will 
have $h^0 (X , \O_X (mD)) >1$ for some $m>0$ (Proposition 4.1 and Lemma 4.3), 
with $m$ explicitly dependent on the class; in particular $mD \sim \Delta + \calE$ for $\Delta $ a movable class  and $\calE$ supported on the rigid non-movable surfaces.  Recall that a class is said to be \it movable \rm if it is in the closure of the cone generated by mobile classes, where an integral class is called \it mobile \rm if it corresponds to a non-empty linear system with no fixed component.  The cone of movable classes is called the \it movable \rm cone $\Movb (X)$.

The main results of this paper are for Picard number $\rho =3$.  A crucial tool in the proofs will the classical Steinian involution on the Hessian of a real elliptic curve.  We recall below that there are only two possibilities where the Hessian of the elliptic curve is singular, and in both these special  cases the real elliptic curve has one real component.

\begin{main} Suppose $X$ is a \CY threefold with Picard number $\rho =3$ for 
which the 
cup-product cubic form and the linear form defined by the second chern class on $H^2 (X, \Z)$ are specified, 
and where the corresponding plane cubic  curve is smooth.  Assume $X$ contains at most one rigid non-movable surface.

(i) If the above real elliptic curve has  two components, then the \Kahler cone of $X$ is contained in the positive 
cone on the bounded component; moreover $X$ lies in a bounded family.

(ii) If the real elliptic curve has one component and the Hessian curve is smooth, then there is a rigid non-movable surface $E$ on $X$ and its class lies in the closed half-cone on which $H \ge 0$ 
determined by the 
bounded component of the Hessian; moreover $X$ lies in a bounded family.

(iii) If the Hessian curve is singular and 
 the line $c_2 =0$ in $\P^2 (\R)$ does not intersect the elliptic curve at an inflexion point, then 
 $X$ contains a rigid non-movable surface $E$ and 
 $X$ lies in a bounded family.
 \end{main}
  
\begin{mainc} If $X$ is a \CY threefold with Picard number $\rho =3$ for which we have specified the above cubic and linear forms on $H^2(X, \Z)$, where 
the corresponding cubic curve and its Hessian  are smooth, then for $b_3(X)$  sufficiently large there must be at least two rigid non-movable surfaces on $X$.
If the Hessian curve is singular, we have the same result provided we assume also 
that the line $c_2 =0$ does not intersect the real elliptic 
curve at an inflexion point.
\end{mainc}

Let me indicate one of the main ideas we shall need.
We will say that a non-zero point $D_0$ on the boundary $\partial W$  of a convex cone $W \subset \R^\rho$ 
is \it visible \rm  from a point $A \in \R^\rho$ if the line segment 
joining $A$ and $D_0$ does not meet the interior of $W$.  If $D_0$ is a smooth point of $\partial W$, the tangent hyperplane 
through $\R_+ D_0$ determines a closed half-space whose intersection with the interior of $W$ is empty;
 then $D_0$ is visible from $A$ if and only if $A$ is in this half-space.  The set of points 
in $\partial W \setminus \{ 0\}$ which are visible from $A$ and $-A$ will be called 
the \it visible extremity \rm of $\partial W$ from $A$; a non-zero smooth point $D_0 \in \partial W$ is in the visible 
extremity from $A\ne 0$ if and only if $A$ is in the tangent hyperplane.
 A non-zero $D_0 \in  \partial W$ which is not in the visible extremity from  $A$ will not be visible from $A$ if and only if it is visible from $-A$, and also if and only if  $-D_0$ is visible from $A$ with respect to $-W$.
A  guiding result is the following:

\begin{prop} Let $X$ be a Calabi--Yau threefold containing finitely many 
rigid non-movable surfaces $E_1 , \ldots , E_r$, with $P^\circ$ is the component of the positive index cone which contains the \Kahler cone; consider the open subcone of $P^\circ$ given by the conditions $E_i \cdot D^2 >0$ for all $i$  and let $Q$ be the component containing the \Kahler cone.
Any point $D_0$ of the boundary of $Q$ at which the cubic form is positive and the Hessian form vanishes, but which is 
 not visible from any of the $E_i$,  must be big.

\begin{proof} We may assume that $X$ is general in moduli, and then the elements of 
the pseodoeffective cone $\Effb (X) $ are those of the form $\Delta + \sum r_iE_i$, for $\Delta$ a real movable divisor and $r_i$ non-negative real numbers.  As commented above, since the \Kahler cone is a subset of $Q$,  any element of $Q$ is of such a form, and so its closure $\bar Q \subset \Effb (X)$.  The claim is that $D_0$ is big.

Suppose to the contrary that $D_0$ lies in the boundary of $\Effb (X)$.  Since it is not visible from any of the classes $E_i$ with respect to $ Q$, it is not visible 
from any of the classes $E_i$ with respect to the convex cone $\Effb (X)$.  
 It follows that the real numbers $r_i$ defined above above must all be zero, and so 
$D_0 \in\Movb (X)$; in fact  $D_0$ must lie in its boundary.  We saw in the second proof of Theorem 0.1 in Section 4 of \cite{WilBd} that for any movable class $D$, we have $\vol (D) \ge D^3$.  Thus $\vol (D_0 ) \ge D_0^3 >0$ and so $D_0$ was big after all. \end{proof} \end{prop}

\begin{rem}
Let us comment how this result will be used.  Assuming $X$ is general in moduli,
we know that $D_0 = \Delta + \calE$, for some real movable divisor $\Delta$ and some real effective class $\calE = \sum r_i E_i$.  In the cases we study, we shall show that  $\Delta \not\in P$ but it is in (the closure of) a different component of 
the positive index cone.  If $L$ denotes an ample class on $X$, then $\Delta + \mu L$ is big and movable for all $\mu > 0$ and hence the Hessian is always non-negative on the line segment from $\Delta$ to $L$; by varying $L$ a little we'll see that the Hessian may be assumed strictly positive on this 
 line segment, apart maybe at $\Delta$ itself.
 As however $(\Delta + \mu L)^3 >0$ for all $\mu >0$, we deduce from  connectedness that $\Delta$ lies in $P$, which will be the required contradiction.
\end{rem}

In particular, assuming there are no rigid non-movable surfaces on $X$, we suppose that $P^\circ$ denotes the 
component of the positive index cone containing the \Kahler cone.
If $X$ is general in moduli, then not only will $Q = P^\circ$ but also any big class will be movable, and in particular has non-negative Hessian.  Thus Proposition 0.1 
implies that  there are  
no points on the boundary of $P$ at which the cubic is positive and at which the Hessian is smooth and vanishes, since otherwise there would be nearby rational points which were both big and had strictly negative Hessian.  
In the case therefore  when when the cubic hypersurface is smooth, this recovers 
the fact that under the above assumptions, the hypersurface must  have two components, and the \Kahler cone on $X$ is contained in the positive (open) cone on the bounded component (\cite{WilBd}, Remark 4.7), which in turn is contained in the \it strictly movable \rm cone $\Mov (X)$, namely the interior of $\Movb (X)$.  Thus any
  integral $D$ in the 
positive cone on the bounded component is big and movable.  

In the case $\rho = 2$, we noted that there are at most two rigid non-movable surfaces, and the more delicate case was 
when there were precisely two --- the case of one such surface was slightly easier.  In this paper we shall mainly be considering the case
for higher Picard number where there is precisely one rigid non-movable surface $E$ on $X$.  In this case we have a strengthening of 
 a result from \cite{WilBd} --- cf. Lemma 5.4 there.  In the case of $\rho =3$ and the cubic is smooth, we deduce that $E\in P$ can only happen if $E$ represents an inflexion point of the elliptic curve.
 
\begin{lem} For arbitrary $\rho >2$,  let  $P$ denote the closure of the component $P^\circ$ of the positive index cone that contains the \Kahler cone and $E$ be the class of the unique rigid non-movable surface on $X$; if $E\in P$, then the cubic form and its Hessian together with the linear form $c_2$ must all vanish at $E$.
 \begin{proof} Given such an $E \in P$, we would have 
  $E^2\cdot D \ge 0$ for all $D\in P$ by Lemma 3.3 of \cite{WilBd}.  
  From the same lemma, we have $E\cdot D^2 >0$ for all $D\in P^\circ$. (If $E\cdot D^2 =0$ for some $D\in P^\circ$, then for any $R$ we have  $E\cdot (D+\varepsilon R)^2 \ge 0$ for $|\varepsilon | \ll 1$ and hence $E\cdot D \cdot R =0$ for any $R$; this in turn implies that $D\cdot E \equiv 0$ and so the Hessian would vanish at $D$.)
So in particular the open subcone $Q$ of  $P^\circ$ defined by $E\cdot D^2 >0$ is the whole of $P^\circ$.
By Proposition 4.4 of the same paper,  the assumptions imply that  $E^3 =0$ and $c_2 \cdot E =0$, since otherwise some multiple of  $E$ would move.  A major part of the proof of 
Lemma 5.4 of \cite{WilBd} was that for any $\rho >1$ there then
 exists a mobile big divisor $L$ with $L\cdot E^2 <0$. 
   In our case therefore,  $E$ represents a smooth point of the cubic form; moreover
    $L$ and $P^\circ$ are on different sides of the tangent hyperplane given by $E^2$.  
  If now the Hessian does not vanish at $E$, then   
  for small $\varepsilon >0$, we have $E - \varepsilon L \in P^\circ =Q$, since $(E -\varepsilon L)^3 >0$ 
 for $0 < \varepsilon \ll 1$.  However, we know then that $E - \varepsilon L = \Delta + \lambda E$ for some movable class $\Delta$ and real 
 $\lambda \ge 0$.  Thus some multiple of $E$ is movable and big; clearly a negative multiple of $E$ cannot be effective, and so we deduce that $E$ is movable, 
 a contradiction.  \end{proof}\end{lem}
 
 Throughout this paper, we may without loss of generality  take $X$ to be general in moduli; recall from \cite{WilKC} that the \Kahler cone is essentially invariant under  deformations and it 
  is always contained in the \Kahler cone of a general deformation.
  
 The overall plan of the paper is as follows.  In Section 1, we describe explicitly the components of the positive index cone associated to a smooth ternary cubic.  We will then consider the case when the cubic comes from the cup-product on $H^2 (X, \Z )$ for $X$ a \CY threefold with Picard number $\rho (X) =3$.  
 For the case when the corresponding real elliptic curve has two components, and $E$ is a rigid non-movable surface on $X$,  we employ the Steinian involution in Section 2 to 
 describe the open subcone defined by the quadratic inequality $E\cdot D^2 >0$ for any component of the positive index cone.  This enables us in this case to deduce 
 that, when $X$ contains at most one rigid non-movable surface,  the \Kahler cone must be contained in the positive cone on the bounded component of the elliptic curve.  We use this in Section 3 to prove boundedness in this case.  Finally in Section 4, we employ the Steinian involution to study in more detail the case when the real elliptic curve has one component.  Under the assumption that there is at most one rigid non-movable surface on $X$, this can occur only for a smooth Hessian when $E$ lies in the closed half-cone on which $H\ge 0$ determined by the bounded component of the Hessian, or 
  in one of the two cases where the Hessian is a singular curve, and in all these cases we can deduce boundedness, in the latter two cases under the mild extra condition that the line $c_2=0$ does not contain an inflexion point of the curve.

\section{Components of the  positive index cone for real ternary cubics}

In \cite{WilBd}, a central role was played by the positive index cone corresponding to the cubic, 
namely the classes $L$ for which $L^3 >0$ and the quadratic form given by $D \mapsto L\cdot D^2$ has index $(1,\rho -1)$.
For $\rho =3$ the cubic defines a curve in the real projective plane, and this is the main case we shall study in this paper.
In this paper, we shall not study the case when the cubic curve  is singular, which would involve case by case arguments,  and our initial aim will be to 
understand the positive index cone in the smooth case.  To study real 
elliptic curves, the Hesse normal form for the curve will be useful, the theory of which may be found in \cite{BM} or Section 3 of \cite{Dolg}.  Normally we might take 
real coordinates so that the real elliptic curve takes  
Hesse normal form  
$x^3 + y^3 + z^3 = 3kxyz$, but for our purposes it will be more convenient to make a change of coordinates so that the cubic takes the form 
$$ F(x,y,z) = -x^3 - y^3  - (z-x-y)^3 + 3kxy(z-x-y), \eqno{(1)}$$
so that the `triangle of reference' of $F$ is now in the affine plane $z=1$ with vertices $ (0,0), (1,0)$ and $ (0,1)$.  We shall write $F_k$ if we wish to indicate the dependence on $k$.  
Recall that if $k>1$, then the real curve $F=0$ has two components, the bounded component (lying in the triangle of reference) and the unbounded component.  The cone 
in $\R^3$ corresponding to the bounded component has two components when one removes the origin, a positive part inside which $F>0$ and a negative part inside which $F<0$, whilst the cone on the unbounded component only has one component in $\R^3$, even after removing the origin.
 In the case $k>1$ it is easily checked that the unbounded component 
 has three affine branches, one of which  
lies in the negative quadrant $x<0,\  y<0$, one in the sector $ y<0,\  x+y >1$ and the third in the sector $x<0,\   x+y >1$.  The inflexion points of the cubic are at 
$B_1 = (0:1:0), B_2 =(1:0:0)$ and $(1:-1:0)$, i.e. the intersection of the line at infinity $z=0$ with the curve (a further reason why the chosen change of coordinates is helpful).  The asymptotes for the affine branches of the unbounded component may be found by calculating the tangents to the curve at the inflexion points, and are  
$$ x = - {1\over {k-1}}, \quad  y = - {1\over {k-1}} \quad \hbox{\rm and}\quad x+y = {k\over {k-1}}. \eqno{(2)}$$

Noting that $x^3 + y^3 + z^3 -3xyz = (x+y+z)(x+\omega y + \omega^2 z)(x+ \omega ^2 y + \omega z)$, where $\omega$ is a primitive cube root of unity, 
 when $k=1$ the cubic (1) splits into the real line $z=0$ and two complex lines 
(meeting at the centroid $({1\over 3}:{1\over 3}:1))$ of the triangle of reference.

When $k <1$, the cubic $F=0$ is smooth but with only one real component, with three affine branches, one in the region $x>0,\ y>0, \ x+y >1$, one in the region 
$x>0, \ y<0, \ x+y <1$ and one in the region $x<0, \ y>0, \ x+y < 1$.  The asymptotes are calculated as before and are given by the equations (2).  The case $k=-2$ is 
special, partly because it is the only smooth case where the asymptotes are concurrent, with the common point being the centroid $({1\over 3}: {1\over 3}:1)$ of the triangle of reference.

If one calculates the Hessian form of the cubic $-x^3 -y^3 -z^3 + 3kxyz$, one obtains 
$$27\bigl(2k^2(x^3 +y^3 +z^3 )- ( 8-2k^3) xyz\bigr).$$
Thus if $H_k$ denotes the Hessian of the cubic $F_k$, the fact that our change of coordinates was unimodular shows for $k\ne 0$ that $H_k = -54k^2F_{k'}$, with constant $k' = {{4-k^3}\over {3k^2}}$.  In particular we see that if $k>1$, then $k' <1$, and so the Hessian curve of a real elliptic curve with two components has only one component.  For $k<1$, we have two notable values: $k=0$ for which the Hessian is the three real 
lines given by $-xy(z-x-y) =0$, and $k=-2$ for which the Hessian curve is the three lines 
(two of them complex) corresponding to $k' =1$ described before.  Apart from these two values, for any real elliptic curve with one component, the Hessian is a real
elliptic curve with two components, the bounded component lying in the triangle of reference.  

The case of the curve having two components is illustrated in Figure 1 (which shows $F_5$, $H_5$ and the
three asymptotes for $F$).  The picture for one component when $k \ne -2, 0$ is not dissimilar, where the roles of the cubic and its Hessian are switched.  
    For $k \ne -2, 0 $ or $1$, we note that the asymptotes to the cubic $F_k =0$ are tangent to the affine Hessian curve $H_k =0$; this is just a special case of a classical result that the 
double polar with respect to the cubic at a point on its Hessian is tangent to the Hessian at the image of the point under the Steinian involution 
(see \cite{Dolg}, Section 3.2 and Exercise 3.8).  When the point is 
an inflexion point of the cubic, the double polar is just the tangent line to the cubic, in our case the asymptote --- the corresponding points on the Hessian are labelled $Q_i$, $i=1,2,3$.  This gives more precise information about the 
affine regions where the Hessian curve can lie.

 \begin{figure}
     \includegraphics[width=10cm]{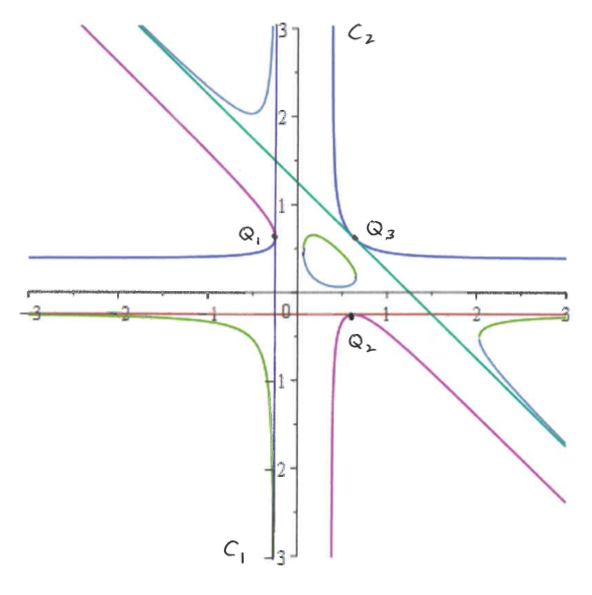}
\caption{Cubic with asymptotes and Hessian for $k=5$}
\end{figure}

We now identify the components of the positive index cone.  Taking the affine piece $z=1$ of a real cubic $F=0$ as above, one looks for the regions for which either $F>0$ and $H>0$, or $F<0$ and $H<0$, the latter being relevant since for $D$ in such a region, both $F$ and $H$ will be positive at $-D$.  In order to keep track of signs, we note that $F_k(0,0,1) =-1$ for all $k$ and that $H_k(0,0,1) >0$ 
for $k\ne 0$.  Moreover the index of the quadratic form $D \mapsto A\cdot D^2$ is $(1,2)$ at $A = (0,0,1)$.

If $F_k =0$ has two real components, i.e. $k>1$, then the curve bounds precisely four (convex) regions of the affine plane $z=1$ on which both $F>0$ and $H>0$.  On these regions, we note that the index 
of the corresponding quadratic form is $(1,2)$.  Moreover the Hessian curve bounds precisely three (convex) regions of the affine plane on which both $F$ and $H$ are 
negative (where the index of the associated quadratic form is $(2,1)$).  Apart from the bounded component inside which $F>0$ (contained in the triangle of reference), each unbounded component 
 defined by $F>0$ will together with the negative of the appropriate region 
defined by $H<0$ give rise to a connected component of the positive index cone in $\R^3$, part of whose boundary is contained in  $F=0$ and part of whose boundary is contained in  $H=0$, with the two parts meeting along rays corresponding to 
 two of the inflexion points of the curve $F=0$.  For each of these three resulting \it hybrid \rm  cones in $\R^3$, we have  $F>0, H>0$ and a continuity argument verifies that the index 
 is $(1,2)$ on each cone.  The other 
component of the positive index cone corresponds to the bounded component and whose boundary is 
contained in $F=0$.

In the case when $F=0$ only has one connected component, i.e. $k<1$, we have that the 
condition $F>0$ defines three (convex) regions of the affine plane $z=1$.  There are two special cases here, when $k=-2$ and $0$.
For all other values of $k<1$, there are four components of the positive index cone.  Three are \it hybrid\rm, obtained from unbounded affine regions on which $F>0, H>0$ together 
with the (negative of) unbounded affine regions on which $F<0, H<0$; a continuity argument again ensures that the index is $(1,2)$.
The remaining component 
corresponds affinely to the bounded component of $H=0$; on the corresponding open half-cone  we have 
$F<0, H<0$, and the index is $(2,1)$ --- this latter claim can be checked in a number of ways (one way is to check it for points in the triangle of reference for the Fermat cubic and using a continuity argument, since as $k\to 0$ the bounded component of the Hessian of $F_k$ tends to this triangle).  Thus points in the negative of the corresponding cone have index $(1,2)$ as required, to give a connected component of the positive index cone.

In the special case $k=0$, the Hessian is $-xy(z-x-y)$; again there will be three \it hybrid \rm components of the positive index cone, each bounded by the cone on one of the affine branches of $F=0$ and two linear parts corresponding to segments of two of the three lines of the Hessian.  As in the general case, there is also the negative cone on the bounded component 
of the Hessian, in this case corresponding to the whole triangle of reference.
The other special case is when $k=-2$: here we have $H>0$ on all the affine plane except at the point $({1\over 3}: {1\over 3}: 1)$, and 
the positive  index cone has just the three connected components, each bounded by the 
(positive)  cone on 
an affine branch of the curve $F=0$ and  a cone in the 
plane $z=0$ corresponding to the line at infinity, the real component of the curve $H=0$.  Here the bounded component 
of $H$ has shrunk to a point, and so does not contribute to the components of the positive index cone.

In both cases, namely the real elliptic curve has one or two components, the components of the positive index cone are all convex, and 
their closures are strictly convex unless $k=0, -2$.

\begin{notation}  We shall now fix on the notation that will be used in the rest of this paper to describe a hybrid component $P^\circ$ of the positive index cone.  
Taking the affine slice $z=1$, where the cubic is of the form described above, with slight abuse of notation concerning the points at infinity, 
we denote the  projectivised boundary of $P$ by $C = C_1 \cup C_2$, where $C_1$ is an 
affine  branch of $F=0$ and $C_2$ is an associated affine branch of $H=0$, 
with $C_1$ and $C_2$ meeting at two inflexion points (at infinity).  When $k>1$, without loss of generality we may take $C_1$ in the negative quadrant and then $C_2$ in the region $x>0, \ y>0, \ x+y >1$. 
The branches $C_1$ and $C_2$ meet (at infinity) at the 
inflexion points $B_1 = (0:1:0)$ and $B_2 = (1:0:0)$.  More specifically, the closed component $P$ then has boundary the positive cone on $C_1$ together with the negative cone on $C_2$, the two parts meeting in  two  rays corresponding to positive multiples of $(0, -1, 0)$ and $(-1,0,0$).  When $k<0$ and $k \ne 0, -2$, we take $C_1$ in the region $x>0,\ y>0,\ x+y >1$ and 
$C_2$ in the negative quadrant. The branches $C_1$ and $C_2$ again meet (at infinity) at the 
inflexion points $B_1$ and $B_2 $.  The component $P$ then has boundary the positive cone on $C_1$ together with the negative cone on $C_2$, the two parts meeting in the two  rays corresponding to positive multiples of $(0, 1, 0)$ and $(1,0,0$).  For $k=0$, the only change is to take $C_2$ to denote the union of the two line 
segments $x\le 0, \ y=0$ and $x=0,\ y\le 0$ in the affine plane $z=1$.  For $k= -2$, we shall take $C_3$ to denote the line segment at infinity from $B_1$ to $B_2$.  The boundary of $P$ in this case is the positive cone on $C_1$ together with the cone generated by $(0,1,0)$ and $(1,0,0)$ in $\R^3$
\end{notation}

\begin{prop} Suppose $X$ is a \CY threefold with Picard number $\rho =3$ and there is at most one 
 rigid non-movable surface $E$ on $X$.  If the cubic form defines an elliptic curve with one real component, for which the Hessian has a non-trivial real  bounded component 
 (i.e. in the notation above, when $k<1$ and $k\ne -2$), then the \Kahler cone of $X$ is not contained in the component  of the positive index cone 
$P^\circ$ corresponding to the bounded component of the Hessian.

\begin{proof}
We suppose the contrary, so that then the convex cone 
 $\Movb (X) \subset P$ by index considerations.  We observed that the result is true if there are no rigid non-movable surfaces on $X$.  
We may assume that $X$ is general in moduli and there is a unique such surface $E$.
We then choose $D_0$ in the boundary of the convex cone $\Movb (X)$ which is not visible from $E$.  Since $\vol (D) \ge D^3$ for all movable classes $D$, we deduce that 
$\vol (D_0) \ge D_0^3 >0$, and hence $D_0$ is big.
We write $D_0  = \Delta + \lambda E$ for some suitable real movable class $\Delta$ and $\lambda \ge 0$, and thus $D_0 - \lambda  E \in \Movb (X)$.  
By convexity  $ D_0 - \mu  E \not \in \Movb (X)$ for all $\mu >0$ since $D_0$ is not visible from $E$, 
hence a contradiction.  \end{proof}\end{prop}

For $X$  a \CY threefold with Picard number $\rho =3$ where there is at most one 
 rigid non-movable surface $E$ on $X$, 
Sections 2 and 4 will rule out the possibility of the \Kahler cone being contained in a \it hybrid \rm component of the positive index cone (except for some special cases 
when the real elliptic curve has one component, for which cases we can prove  boundedness results), 
where \it hybrid \rm by definition means that part of the boundary consists of points where $F$ vanishes and part where $H$ vanishes, as described above.  This then reduces us down in Section 3 to 
considering  the case where the elliptic curve has two components and the \Kahler cone is contained in the positive cone on its bounded component,  for which case we can 
also prove boundedness  (Corollary 3.8).  

We note that if $D_0$ is a point in the open arc of rays in the boundary of $Q$ as specified in the statement of the following Proposition,  then  a line segment from  $E$ to the 
given class $D_0$ contains interior points of $Q$ if and only if it contains interior points of $P^\circ$.  
 So in this case we can 
talk about $D_0$ being visible or not visible from $E$ without specifying whether we consider $D_0$ to be on the boundary of $Q$ or $P^\circ$, or 
indeed say the negative cone on the interior of $C_2$.

\begin{prop} Suppose $X$ is a \CY threefold with Picard number $\rho = 3$ for which the cubic form defines a smooth elliptic curve, and that $X$ contains a unique rigid non-movable surface $E$.  Let $P^\circ$ denote a component of the positive index cone and $Q$ a component of 
the open subcone of $P^\circ$ given by $E\cdot D^2 >0$.  If the boundary of $Q$ contains a non-empty open 
arc of rays which are not visible from $E$ but on which the Hessian 
$H$ vanishes and the cubic form $F$ is strictly positive, 
 then $Q$  does not contain the \Kahler cone.

\begin{proof} As usual, we may assume that $X$ is general in moduli.  In the light of Proposition 1.1, we may assume $P^\circ$ to be a hybrid component of the positive index cone; we assume that $Q$ contains the \Kahler cone and 
obtain a contradiction.
Let $D_0$ be any point in  one of the specified rays in the boundary;   Proposition 0.1 then implies that for some $\lambda \ge 0$, the class $D_0 - \lambda E$ is movable.  Moreover, the argument from the proof of Proposition 1.1 shows that $D_0$ is not movable and $\lambda >0$.  
With the notation from this section, we may
take $-D_0 \in C_2$.  Let us assume first that  we are not in either of the two special cases $k=0$ or $-2$ --- the argument will in fact be very similar in these cases.

We then have the possibilities that either some real multiple of $E$ represents a point $A$ of the affine plane $z=1$, or $E$ represents a point 
 on the line at infinity $z=0$.  Suppose first that $A$ is a negative multiple of $E$, and so $H(A) \le 0$.  Then the point $-D_0 \in C_2$ is not visible with respect to the cone  $-Q$ from the point $A$, and in particular  the line segment from $-D_0$ to $A$ 
cuts the Hessian curve at points $-D_0$, $A' \in C_2$ and one other point; note here that we are using $E\not \in P$ from Lemma 0.3.  Consider now $\Delta = -\lambda E +D_0$; if this were to lie in the half-plane 
$\{ z\ge 0\}$, then we would have $H(\Delta ) <0$, a contradiction.  Thus $\Delta$ lies in the half-plane $\{z\le 0 \}$, and so $H(-\Delta ) \le 0$, showing that $\Delta$ is in the closure of a different component of the positive index cone.

The next possibility is that $A$ is a positive multiple of $E$; then $-D_0 \in C_2$ is visible from $A$; moreover 
by Proposition 0.1 we have $\lambda >0$ such that $\Delta = D_0 -\lambda E$ is movable, and hence for some $\mu >0$ that $\Delta = D_0 -\mu A$ is movable.  Now consider the line segment consisting of points $-D_0 + t A$ for $0\le t \le \mu$; since $H(-D_0 + \mu A) \le 0$, the Hessian will vanish when $t=0$ and precisely one more point.  We deduce therefore that at $-\Delta$ we have $H\le 0$ and $F<0$, so that $\Delta$ is in 
(the closure of) a different component of the positive index cone.  

Finally, if $E$ represents a point on the line at infinity, we know from Lemma 0.3 
that either $E$ represents one of the two inflexion points at infinity on $C_1$ 
(noting that the sign is determined since $E\cdot H^2 >0$ for any ample class $H$), or 
$E\not \in P$.  In the first case, since $-\Delta = -D_0 +\lambda E$ and $H(-\Delta) \le 0$, we see from the geometry of the affine curves that 
$\Delta$ is in 
(the closure of) a different component of the positive index cone.  If $E$ represents the third inflexion point (and so $E\not\in P$), then we must have $Q= P^\circ$; the fact that $-D_0$ is visible from $E$ implies (again from the affine geometry) that $H(-D_0 +tE) > 0$ for all $t>0$, contradicting the fact it is non-positive for $t=\lambda$.  The remaining possibility is that $E$ lies in the interior of the cones on the other two line segments with $z=0$ on which $F>0$ (since we know that $H(E)\ge 0$).  The fact that the chosen class $-D_0$ is visible from $E$ implies that it is visible from $E$ with respect to the cone on $C_2$, and the affine geometry then again yields that 
$H(-D_0 +tE) > 0$ for all $t>0$, contradicting the fact it is non-positive for $t=\lambda$.

In all possible cases therefore, $\Delta$ is in 
(the closure of) a different component of the positive index cone, and so in particular $\Delta ^3 \ge 0$, 
and we 
obtain a contradiction via the argument in Remark 0.2.

The reader is invited to check that the same arguments hold true in the case when $k=0$ as for the general case when $k<1$;  the ability to perturb $D_0$ ensures that in the case where the  movable divisor $D_0 - \lambda E$ is in the closure of another component  of the positive index cone, it is not also in $P$.    For $k=-2$, it is even more straightforward;   in this  case 
 $H(D) >0$ if and only if $D$ is a point of the affine plane $z=1$ and not the singular point.  We know also that the open arc specified in the Proposition corresponds to an open interval 
 of the line segment $C_3$ at infinity, i.e all its points lie in the cone generated by $(0,1,0)$ and $(1,0,0)$.  We have seen in Lemma 0.2 that $E\not\in P$ unless it is one of the points at infinity on $C_1$, and since $E\cdot L^2 >0$ for any ample $L\in P$, it follows from Lemma 3.3 of \cite{WilBd} that $-E \not \in P$.  Moreover,
 if $E$ were to lie on the line at infinity, then any choice of $D_0$ in the open arc specified would be visible from $E$.
 The assumptions of the Proposition therefore ensure that $E$ is a positive multiple of a point $A$ in the affine plane $z=1$, and given $E$ we may perturb our choice of $D_0$ if required so that 
 the Hessian $H$ is non-singular at all points on 
 the line through $D_0$ and $-A$; hence 
the Hessian  takes negative values on $D_0 - 
\lambda E$ for all $\lambda >0$;  thus $D_0 -\lambda E$ is not  movable for any $\lambda >0$, contradicting the conclusion of Proposition 0.1.  \end{proof}\end{prop}

\section{Steinian map when elliptic curve has two real components}

With notation as in the previous section, we let  $P^\circ$ denote a component of the positive index cone, with closure $P$.
The next piece of the jigsaw puzzle will be to understand the function $G_A (x,y) = A\cdot D^2$, where $A = (a,b,1)$ and $D= (x:y:1)$, and in particular where it vanishes on 
the affine slice of the boundary of a component  $P$ of the positive index cone.  
In this section, we study the case when $k>1$, i.e. the real elliptic curve $F=0$ has two components.

The component of the positive index cone 
which corresponds to the bounded component of $F=0$ is then easy to understand; phrasing things affinely,   
we know that $A\cdot D^2 =0$ for $D$ on the bounded component is just saying that $A$ lies on the tangent to the curve at $D$.  
 Assuming that $A$  corresponds to a point 
outside the bounded component, the conic $G_A =0$ 
intersects that boundary of the component precisely at the points where the tangent to the curve passes through $A$, and hence 
corresponds to the boundary of the \it visible points \rm from $A$ 
(where we recall \it visible \rm here means the points $D$  on the bounded component for which the closed line segment $AD$ does not include any points inside the component); we used the term 
 \it visible extremity \rm for the set of such points.

We therefore need to understand the hybrid components $P^\circ$.  Here we adopt the notation from Section 1.  In view of Lemma 0.3, we assume $A\not\in P$; 
  as above there is then a simple answer to the question 
of where on  $C_1$ we have vanishing of $G_A$, namely points on of $C_1$ for which the tangent passes through $A$ (including maybe inflexion points at infinity).  There will be two such points if $A$ is in the quadrant 
$x\le -{1\over {k-1}},\  y\le -{1\over {k-1}}$ (including the possibilities of the inflexion points $B_1$ and $B_2$), 
no such points if $A$ is in the quadrant 
$x> -{1\over {k-1}}, \ y> -{1\over {k-1}}$ and one point otherwise.   Let us consider the branch $C_2$ of the Hessian $H$ passing through $B_1, B_2$ and $R = Q_3$, the affine point 
with coordinates $( {k\over {2(k-1)}}, {k\over {2(k-1)}})$.  As $B_3 = (-1:1:0)$ is the third inflexion point, we have already observed that the tangent to $F=0$ at $B_3$ is  tangent to $H=0$ 
at $R$, i.e. that if we take $B_3$ to be the zero of the group law, then $R$ is 
the unique real 2-torsion point of the Hessian.  Recall the classical fact that there is a 
well-defined base-point free involution $\alpha$ on the Hessian 
curve, known as the \it Steinian involution \rm (\cite{Dolg}, Section 3.2, noting a misprint in Corollary 3.2.5), where the polar conic of $F$ with respect to a point $U$ on $H$ is a line pair with singularity at
$\alpha (U)$.  In the case currently under consideration where the Hessian has only one real component,  
corresponding to a choice of inflexion point for the zero of the group law, there is a  
unique real 2-torsion point, and $\alpha$ is given by translation in the group law by this point.  The Steinian map has the property 
that for any point $U \in H$, the second polar of $F$  with respect to 
the point $U$ is the tangent to the Hessian $H$ at $\alpha (U)$ (\cite{Dolg}, Exercise 3.8, noting a misprint).  Let $Q_1$ be the point on the branch of $H$ in the region $x<0, \ y>0, \ x+y <1$, which is the 2-torsion point when we take $B_1$ as the zero in the group law,  $Q_2$ the point on the branch of $H$ in the region $x>0, \ y<0, \ x+y <1$ 
corresponding to $B_2$, and $R = Q_3$ the point on the branch of $H$ in the region $x>0,\  y>0, \ x+y >1$ 
corresponding to $B_3$.  It is left to the reader to check that $\alpha (B_1) = Q_1$, $\alpha (B_2) =Q_2$ and $\alpha (R) = B_3$.  Thus the second polar of $F$ with 
respect to the inflexion point $B_1$ is the tangent to the Hessian at $Q_1$, namely given affinely by $x= -{1\over {k-1}}$, the 
second polar of $F$ with 
respect to the inflexion point $B_2$ is the tangent to the Hessian at $Q_2$, namely given affinely by $y= - {1\over {k-1}}$,  and the second polar of $F$ with 
respect to the point $R$ is the tangent to the Hessian at $B_3$, namely the asymptote $x+y = {k'\over {(k'-1)}}$, where $k' = {{4-k^3}\over {3k^2}}$ as defined before.

Setting $A =(a,b,1)$, it is easily checked from this that $G_A (B_1) >0 $ iff $a< -{1\over {k-1}}$, that $G_A (B_2) >0$ iff $b < -{1\over {k-1}}$ and 
that $G_A  (R)>0$ iff $a+b < {{k'}\over {(k'-1)}}$.  Moreover if $A=Q_i$, then $G_A (B_i) =0$ for $i=1,2$ and if $A= B_3$ then $G_A(R)=0$. We are interested in the cases of $B_1$ and $B_2$ 
since we want to know the sign of $G_A$ on points of $C_2$ where either $y \gg 0$ or $x\gg 0$.  
This 
gives us a dictionary as to how many points of $C_2$ there are at which  the function $G_A$ vanishes.  

\begin{lem} Suppose $a< -{1\over {k-1}}$ and $b > -{1\over {k-1}}$.  If 
$a+b > {{k'}\over {(k'-1)}}$, so that $G_A(R) <0$, then $G_A$ vanishes 
precisely once on the interior of the open arc $B_1 R$ of $C_2$, and does not vanish on the interior of the open arc $RB_2$.
If $a+b  = {{k'}\over {(k'-1)}}$, then $G_A$ vanishes once on the interior of the arc $B_1 R$ and twice at $R$, and is non-vanishing on the open arc $RB_2$.
 If $a+b < {{k'}\over {(k'-1)}}$, so that 
$G_A(R)>0$, and $A$ lies in the open region between the arc of the Hessian from $Q_1$ to 
$B_3$, the line $x = -1(k-1)$ and the line $x+y = k'/(k'-1)$, then $G_A$ will vanish precisely twice on the 
open arc $B_1R$ and once on the open 
arc $RB_2$.  If $A$ lies on the arc of the Hessian from $Q_1$ to $B_3$, then $G_A =0$ is a line pair whose 
singularity $D = \alpha (A)$ lies in the interior of the arc $B_1 R$ of $C_2$, one line of the pair also intersecting $C_1$ and the other line also intersecting $C_2$ at some point in the interior of the arc $R B_2$.
For the remaining points $A$ 
satisfying the initial inequalities, $G_A$ is
non-vanishing on the arc $B_1R$ and has one zero on the interior of the arc $RB_2$ of $C_2$.

\begin{proof} 
Note that the second polar of a point $D$ in the arc 
of $C_2$ from $B_1$ to $R$ corresponds to the tangent line to the Hessian at $D' =\alpha (D)$ on the arc of the Hessian going from $Q_1$ to $B_3$.  Furthermore we note that for a given point $A$, we are asking how many points $D'$ on the arc of the Hessian from  $Q_1$ to $ B_3$, the tangent at 
$D'$ contains $A$.  We note 
that $A$ is  
clearly on no such tangent line if $A$ is strictly below the  arc of the Hessian going from $Q_1$ to $B_3$, and if $A$ lies on this arc, the answer is one 
(as the relevant tangent line is just the tangent to the Hessian at $A$).  
In this case $G_A =0$ is a line pair with singularity at $\alpha (A)$ as claimed.   From 
the geometry we see that $A$ also lies on a tangent line at a point of the Hessian arc $B_3Q_2$ and also 
on a tangent line of $C_1$; this then gives the required statements in this case.
If $A$ is strictly above the arc $Q_1B_3$, then the claims are clear 
geometrically since we are looking for 
 points of the visible extremity of the convex set with boundary the affine branch $B_2B_3$ of the Hessian
as seen from $A$, lying on the arc of the Hessian from $Q_1$ to $B_3$, and there will be one such point 
unless $A$ lies in the open region specified in the statement of the Lemma, in which case there are two.  
The remaining claims of the lemma are also clear from the geometry, since we are additionally asking 
if there are points on the open arc $B_3Q_2$ of the Hessian, the tangent line at which passes through $A$.\end{proof}\end{lem}

By symmetry, we have a corresponding result for  $a> -{1\over {k-1}}$ and $b < -{1\over {k-1}}$.  We can then deduce the following result.

\begin{cor} If $a \le -{1\over {k-1}}$ and $b \le -{1\over {k-1}}$, then $G_A$ vanishes twice on the closed arc $C_1$ (including possibly the inflexion points $B_1, 
B_2$ in the case of an equality),
and not at all on the affine branch $C_2$.  If $a > -{1\over {k-1}}$ and $b > -{1\over {k-1}}$ then $G_A$ is non-zero on the closed arc $C_1$;  in this case if  $a+b > {k'\over {(k'-1)}}$, then
$G_A$ will also be non-zero on the arc $C_2$, and if $a+b \le {k'\over {(k'-1)}}$ it will vanish twice on $C_2$ (in the case of equality, twice at the point $R$).  In all other cases, 
$G_A$ will intersect (the closure of) $C_1$ precisely once and the affine arc $C_2$ either once or three times 
\end{cor}
Summing up therefore, the conic $A\cdot D^2 =0$ either intersects the projectivised boundary of $P$ twice (including the case when it is the point $R$ with multiplicity two), or four times in the cases specified. 
 Since we know that $P^\circ$ is convex, we know in particular that 
the  cone  $P^\circ \cap \{ D \  :\  A\cdot D^2 >0 \}$ is connected when we just have the two points of intersection, but will have two components when there are four (distinct) intersection points.
If the index at $A$ is $(1,q)$ with $q\le 2$, 
it follows from Proposition 3.4 of \cite{WilBd} that  any component $Q$ of the above cone 
 is also convex. 
If $A$ has index $(2,1)$, then $-A$ has index $(1,2)$, and the cone 
$Q_- =P^\circ \cap \{ D \  :\  A\cdot D^2 <0 \}$ is convex.  
If however $A$ lies on the affine arc of the Hessian from $Q_1$ to $B_3$, the second polar with respect to $F$ is a line pair with a
singularity at the corresponding point $\alpha (A)$ of $C_2$, as detailed in Lemma 2.1.  The open cone 
$P^\circ \cap \{D: D^2\cdot A >0 \}$ then has two (convex) components, and the open cone 
$P^\circ \cap \{D: D^2\cdot A <0 \}$ has one (convex) component, reflecting the fact that if the index at $A$ is $(1,1)$, then that is also the index at $-A$.

\begin{prop} Suppose $A = (a, b,1)$ is in the affine plane $z=1$ with $a < -{1\over {k-1}}$, $b > -{1\over {k-1}}$. If $H(A) \ge 0$, then any component of $C_2 \cap \{ A\cdot D^2  >0 \}$ contains an open arc of points visible 
from $A$.  
In the case of $H(A) \le 0$, some open arc in $C_2 \cap \{ A\cdot D^2  <0 \}$ is not visible from $A$.

\begin{proof} Since $G_A(B_1)>0$, we have that $G_A = 0$ intersects $C_2$ precisely once on the open arc $B_1R$ if $a+b \ge {k'\over {(k'-1)}}$, and in this case $H(A) >0$.  The condition $A\cdot D^2 \ge 0$ defines a convex subcone of $P$ containing $(0,-1,0)$, part of whose boundary 
corresesponds to an arc of $C_2$ in $B_1 R$ (namely the part near $B_1$), which is 
 visible from $A$.  If  $a+b < {k'\over {(k'-1)}}$, then we may have $H(A) \ge 0$ or $H(A)\le 0$.  
 If $H(A)\ge 0$, then the condition $A\cdot D^2 \ge 0$ defines a subcone of $P$ containing $(0,-1,0)$, and 
 part of the boundary of one component 
corresponds to an arc in $C_2$, part of which (namely the part near $B_1$) is 
 visible from $A$. If the above subcone contains a second component, then part of its boundary corresponds to a second arc in $C_2$ whose interior contains $R$, and so part of that arc is also visible from $A$.
 If $H(A) \le 0$, then the condition $A\cdot D^2 \le 0$ defines a convex subcone of $P$ containing $(-1,0,0)$,
 part of whose boundary corresponds to an arc of $C_2$ in $RB_2$, part of which (namely the part near $B_2$)
is not visible from $A$.\end{proof}\end{prop}

 \smallskip
 
 \begin{figure}
 \includegraphics[width=9cm]{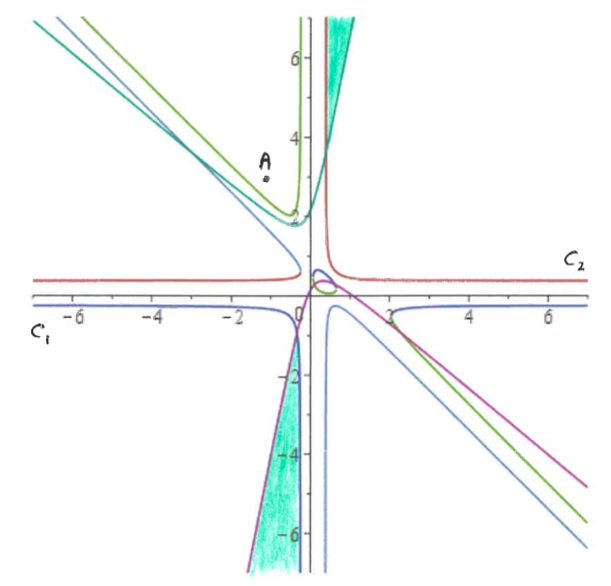}
\caption{Case $k=5$ and $A=(-1,3,1)$}
\end{figure}

Before proceeding further, the reader might find Figures 2 and 3 helpful.  
      We are taking $k=5$ as in Figure 1, and $C_2$ is the branch of the Hessian in the positive quadrant.  In both figures, the point $A$ is marked, 
      $P^\circ$ is the hybrid component of the positive index cone with boundary corresponding to $C_1 \cup C_2$, 
      and the areas corresponding to the subcone $Q$ have been shaded.
For Figure 2 we take $A$ to be the affine point $(-1, 3, 1)$; here $H(A) >0$ 
and so by Proposition 2.3, we have part of $C_2 \cap \{ A\cdot D^2  >0 \}$ is visible from $A$.   The picture 
when $A$ lies in the open region between the arc $Q_1B_3$ of the Hessian and its asymptote is similar to 
Figure 2, except that $P^\circ \cap \{ A\cdot D^2  >0 \}$ has two possible components, with one having part of its boundary corresponding to an open arc of $C_2$ containing $R$, part of which is visible from $A$.
For Figure 3 we take  $A= (-2,1,1)$, so that  
$H(A) <0$ and part of $C_2 \cap \{ A\cdot D^2  >0 \}$ is not be visible from $A$.

\smallskip 
 \begin{figure}
      \includegraphics[width=9cm]{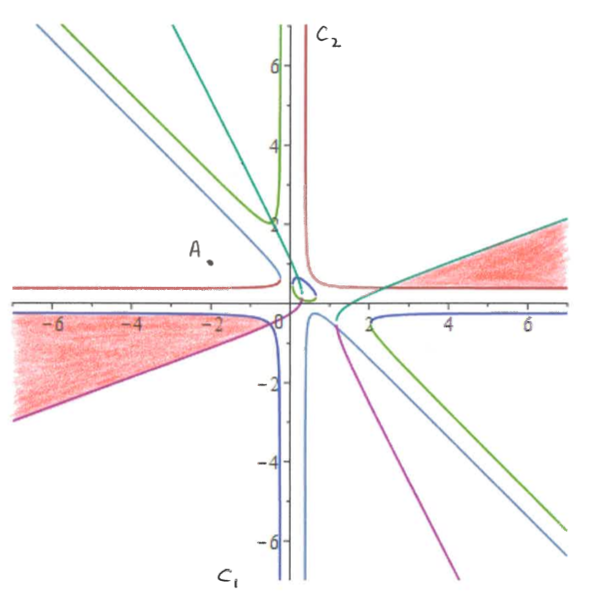}
      \caption{Case $k=5$ and $A=(-2,1,1)$}
      \end{figure}
    \smallskip

\begin{prop} If $X$ is a \CY threefold with Picard number $\rho =3$ whose corresponding real cubic curve is smooth with two real components, and there is at most 
one rigid non-movable surface $E$ on $X$, then the \Kahler cone must be contained in the positive cone on the bounded component of the elliptic curve.

\begin{proof} The content of this result is that the \Kahler cone cannot be contained in a hybrid component of the positive index cone.
We already know this in the case of no such surfaces, so suppose we have a unique such surface $E$. We suppose that the result is not true, so that the 
\Kahler cone is contained in a hybrid component $P^\circ$, which we may assume to be bounded by $C_1$ and $C_2$, with the notation as before.  The index of $E$ at the point $E$ is $(1,q)$ for 
$q\le 2$.  We let $Q$ denote the convex cone $P^\circ \cap \{ E\cdot D^2 >0 \}$.

\begin{claim} There is an arc of points on $C_2$, the negative multiples 
of which correspond to rays as specified in the statement of Proposition 1.2.
\end{claim}

Assuming $E$ does not lie on the line at infinity $z=0$, we 
 let $A$ denote the intersection of the line though $E$ with the affine plane $z=1$.  
Suppose the \Kahler cone is contained in a hybrid $P^\circ$, which without loss of generality we may take to comprise of 
the (positive) cone on interior of the affine branch $C_1$ of $F=0$ together 
with the negative cone on the interior of the affine branch $C_2$ of $H=0$, the two parts  meeting along the cone spanned by $(0,-1,0)$ and $(-1,0,0)$.

If $A$ is such that $G_A =0$ does not properly intersect either $C_1$ or $C_2$, namely $A=(a,b,1)$ with $a> -1/(k-1), \ b> -1 /(k-1)$ and $a+b \le k'/(k'-1)$, 
then $G_A$ is negative on the whole cone $P^\circ$, which contradicts $P^\circ$ containing the 
\Kahler cone, since in this case we note that the index of $A$ is $(1,2)$ and so $A$ is a positive multiple of $E$.
If $A$ is such that $G_A =0$ intersects $C_1$ twice (possibly including points at infinity), then we check easily 
that the boundary of $Q$ contains the whole negative cone on $C_2$.  In this case we also note that the index of $A$ is $(1,2)$ and so $A$ is a positive multiple of $E$.  Since
all of  $C_2$ is visible from the affine point $A$, none of the negative cone on $C_2$ may be seen from $E$, and so the corresponding points of the boundary of $Q$ may not be seen from $E$, and the Claim is proved.
  
If now $G_A$ vanishes twice on $C_2$ (namely $a> -1/(k-1), \ b> -1 /(k-1)$ and $a+b \ge k'/(k'-1)$), then $G_A$ (which is negative at $B_1$ and $B_2$) will be positive on a region inside of $C_2$ with boundary including 
the arc of $C_2$ between the two intersection points.  Such a point $A$ is also a positive multiple of $E$.
Points of this arc may plainly be seen from $A$, and so the corresponding points of the boundary of $Q$ may not be seen from $E$, and the Claim follows.

We are left with the case where $A$ is as in Proposition 2.3, or by symmetry where the roles of $a$ and $b$ are switched.  The Claim then follows from this previous result.

Finally we should mention the case when $E$ lies on the line at infinity.  
We know that $E \not \in P$ by Lemma 0.3, unless $E$ represents one of the inflexion points $B_1$ or $B_2$.  In this 
latter case $E\cdot D^2 >0$ on $P^\circ$, and so there is no difficulty finding an arc of rays in the boundary 
of $P$ which are not visible from $E$.
 If $E$ represents the inflexion point $B_3$, then $E\cdot D^2= 0$ touches the boundary of $P$ along two rays, corresponding to $R$ and the symmetry point on $C_1$.  Thus the subcone $Q$ is either empty or the whole of $P^\circ$, in which latter  case the claim is true since half of the cone 
on $C_2$ is visible from $(-1,1,0)$ and half is not.    The other case is when $H(E)>0$; knowing that $E\not \in P$, we have that 
$E = (a,b,0)$, where \it either $a<0, a+b >0 $ or \rm $b<0, a+b >0$, where by symmetry we may assume the former case.  In this case, the same argument says that $E\cdot D^2$ is 
positive on an arc of $C_2$ (namely from $B_1$ to the point on the arc $B_1 R$ where $G_E$ vanishes), which points are visible from $E$ and again the Claim follows.

Thus we merely need to combine the Claim with Proposition 1.2 to deduce that the \Kahler cone cannot be contained in a hybrid component of the positive index cone, 
and so it must be contained in the  positive cone on  the bounded component of the elliptic curve. \end{proof}\end{prop}

\section{Boundedness when elliptic curve has two real components}

   The result of Proposition 2.4  extends to give a corresponding result for any $\rho \ge 3$.  Note however that the cubic having two real  components in higher dimensions is a stronger condition; for smooth cubic surfaces
for instance, it is well-known that there 
are five real diffeomorphism  types and only one of these has two components.

\begin{thm}  Let $X$ be a \CY threefold which is general in moduli with Picard number $\rho \ge 3$ whose corresponding cubic hypersurface is smooth with two real components.  If there is at most 
one rigid non-movable surface on $X$, then the \Kahler cone must be contained in the positive cone on the bounded component of the cubic hypersurface.

\begin{proof} We choose a general ample integral class $L$ and a general integral class $J$ inside the positive cone on the bounded component of the cubic hypersurface, 
and we take the 3-dimensional slice $\Lambda$ of $H^2(X, \R ) $ through the origin and the classes $E$, $L$ and $J$.  We therefore obtain a cone on a smooth elliptic curve 
with two real components, defined by the restriction of the cubic form to $\Lambda$.  Unless $\rho =3$, the Hessian of the elliptic curve is not of course the restriction of the Hessian hypersurface to $\Lambda$ (if only since the degrees are different).  If $D\in \Lambda$, then the quadratic form on $\Lambda$ associated with a non-trivial divisor $D$ has 
index $(1,2), (2,1)$ or $(1,1)$, according to whether the Hessian of the (smooth) restricted cubic is positive, negative or vanishes.  It will turn out that the Hessian of the restricted cubic plays a more crucial role than the restriction of the original Hessian.

We suppose that the \Kahler cone $\calK$ is not contained in the positive cone on the bounded component, and so $\calK \cap \Lambda$ is not contained in the positive cone on the bounded component of the elliptic curve; we suppose that $\calK \cap \Lambda$ is contained in some other component $P^\circ$ of the positive index cone of the elliptic curve   
--- the restriction of the component of the positive index cone of the cubic on $H^2 (X, \R )$ containing $\calK$ is an open subcone of $P^\circ$. 
We let $Q \subset P^\circ \subset \Lambda$ denote a component of the open subcone of $P^\circ$ given by the extra condition $E\cdot D^2 >0$, which we saw in the previous section was convex.  The proofs of Proposition 4.1 and Lemma 4.3 
from \cite{WilBd} extend to give that 
for any choice of integral $D$ in $Q$, we can find an explicit positive integer $m$ depending on $D$ for which $h^0 (X, \O _X (mD)) >1$, i.e $|mD| = |\Delta | + rE$ with 
$r\ge 0$ and $\Delta \in 
\Lambda$ having some multiple mobile.  Any irreducible surface class in $\Lambda$ will still define a quadratic form with index $(1,q)$ for $q=1$ or $2$ with respect to the cubic form restricted to $\Lambda$;
in particular this holds both 
for $E$ and any non-trivial mobile class in  $\Lambda$, and so the Hessian of the elliptic curve is non-negative there.

The argument of Proposition 0.1 and the results of the previous section applied to the restriction of the cubic to $\Lambda$ then go over directly to yield a contradiction. \end{proof}\end{thm}

Let us consider the case when the hypersurface corresponding to the cubic has two components, and that the linear form $c_2$ has also been specified.  Given this 
information, and assuming that $X$ is general in moduli and contains only one rigid non-movable surface, we aim to prove boundedness.  

We now let $P^\circ$ denote the component of the positive index cone which contains the \Kahler cone.
We have seen in Theorem 3.1 that 
 $P^\circ$ must be  the positive cone on the bounded component of the cubic hypersurface.  The Hessian does not vanish on $P$ (apart from trivially at the origin) and from the case of elliptic curves, we know that $P$ is a strictly convex cone.   
 If $E$ denotes the unique rigid non-movable surface on $X$, it does not represent 
 a point in $P$ by Lemma 0.3.  
 Therefore, away from the origin,    
 the quadric $E\cdot D^2 =0$ intersects $P$ precisely on the visible extremity from $E$, the lines from $E$ to such points being tangent to $P$.  
 This enables us to produce a version of the $\rho =2$ result from Proposition 5.1 in \cite{WilBd}, which follows since $P^\circ$ is contained in the cone generated by $Q$ and $E$, 
 where $Q = P^\circ \cap \{ D\ :\ E\cdot D^2 >0 \}$.

\begin{cor} With the above assumptions, namely $X$
 general in moduli  with $\rho \ge 3$ and with one rigid non-movable surface on $X$, where the cubic corresponds to a smooth 
hypersurface with two real components 
and $P^\circ$ denotes the component of the positive index cone which contains the \Kahler cone (by Theorem 3.1 the positive cone on the bounded component), 
 it then follows that $P\subset \Effb (X)$. \end{cor}

If we know $E$, we choose any integral $D = \Delta +\lambda E$ in the interior of $Q$ with some multiple of $\Delta$ mobile and then bound $\lambda$ (cf. 
\cite{WilBd}, Proposition 5.5) to get only finitely many possibilities for the class of $\Delta$.

\begin{prop} Under the assumptions of the previous result, there is an upper  bound $\lambda _0$ depending only on the known information, including the class of $E$ and the choice of integral divisor $D \in Q$, on the possible $\lambda$ for which $D-\lambda E$ is movable.  There is an integer $m>0$ depending only on $D$ 
and other known data such that 
$h^0 (X, \O_X (mD))>1$, so that  $| mD| = |\Delta | + r E$, where some multiple of $\Delta$ is mobile and $r$ has a known upper bound.  By appropriate choices of $D$ and $m$, we can find a finite set of classes, at least one of which will be movable and big.

\begin{proof} Let us assume first that 
$E^3 \ge 0$.  Since  $E\not\in P$ by Lemma 0.3, it follows that for any $L\in Q$, there exists $\mu = \mu (L) >0$ such that 
$$ 0 = (L+\mu E)^3 = L^3 + 3\mu L^2\cdot E + 3 \mu^2 L\cdot E^2 + \mu ^3 E^3.$$
As three of the terms on the righthand side are positive, we deduce that $L\cdot E^2 <0$ for all $L\in Q$.  This argument also works for $L \in \bar Q$, except for the points on the visible extremity of $P$ from $E$.  Suppose that $L$ is one of these points, so that $L^2 \cdot E =0$ and $\mu =0$.  
Since $L\cdot E^2 \le 0$ by continuity, let us assume that $L\cdot E^2 =0$.  Using the Hodge Index theorem on $E$, we have equality in $E^3 (L^2\cdot E) \le (L\cdot E^2)^2 = 0$, and so $L|_E 
\equiv 0$; hence the 1-cycle $E\cdot L$ in numerically trivial.
This contradicts the fact that the Hessian does not vanish at $L$.  Thus we deduce that in fact $L\cdot E^2 <0$ for all $L\in \bar Q$, i.e. the linear form 
$E^2$ is strictly negative on $\bar Q$.

Having chosen an integral $D \in Q$, there exists $\lambda _0 >0$ depending on the given data (including $E$ and $D$) such that $(D-\lambda E)^2 = D^2 -2\lambda D\cdot E + \lambda ^2 E^2$ is strictly negative on $\bar Q$ for $\lambda \ge \lambda _0$.  In particular, if $L$ is an ample class in $Q$, we have $(D-\lambda E)^2\cdot L <0$, meaning that $D-\lambda E$ cannot be movable; so when $E^3 \ge 0$, we have found a suitable upper-bound 
$\lambda _0$ depending on known data.

Assume now that $E^3 <0$, then plainly $(-E)^3 >0$ and 
the cubic $(D-tE)^3 $ is  positive for $0< t< \lambda _1$, negative for $\lambda _1 <t< \lambda _2$, 
and positive again for $t> \lambda _2$, for suitable positive numbers $\lambda _1 < \lambda _2$, where $t=\lambda _1$ corresponds to leaving the positive cone on the bounded component.  
Thus if $\Delta = D-\lambda E$ is movable, either we have $\lambda \le \lambda _1$, or 
$\lambda >\lambda_2$ and  the mobile class $\Delta$ has $\Delta ^3 >0$.  In this latter case however, if we 
  let $L$ denote an ample class in $P$, then clearly  
$(\Delta + tL)^3 >0$ for all $t\ge 0$, from which we deduce that $\Delta \in P$, a contradiction since then 
$\lambda \le \lambda _1$.  
 Thus in this case too  we have an upper bound $\lambda _2$ on the $\lambda$ for which $D-\lambda E$ is movable.

Now recall from the proofs of Proposition 4.1 and Lemma 4.3 of \cite{WilBd} that there exists $m_0 >0$, 
depending only on $D$ and the cubic and linear forms on $H^2 (X, \Z)$, such that $h^0(X, \O_X (mD))>1$ for $m\ge m_0$, and so $|mD| = |\Delta | + rE$, with some multiple of $\Delta$ 
mobile (and hence $\Delta$ movable), and $r$ a non-negative integer. Therefore we deduce bounds 
for $r$ that $r< m \lambda _0$ if $E^3\ge 0$ and $r< m \lambda _2$ if $E^3< 0$.  Thus there exist a finite set of classes, at least one of which will be a movable 
class $\Delta$, some multiple of which is mobile.    

The movable class $\Delta$ obtained in this way corresponds on some minimal model to a morphism; if $\Delta$ is not big, then the image of this morphism is a surface or a curve.  We can then however choose another integral class $D' \in Q$, not representing a point on the plane through the origin containing $E$ and $D$.  Repeating the argument, we achieve a finite set of possibilities for a second 
movable class $\Delta '$, not a multiple of $\Delta$, and even if both $\Delta$ and $\Delta '$  correspond to maps 
to curves, their sum would correspond to a map to a surface. We may assume then that the movable 
class $\Delta = mD - rE$ obtained is either big or corresponds to a map to a surface.  In the latter case, 
some minimal model of $X$ is an elliptic fibre space over a surface $S$, and the real Cartier divisors on $S$ have transforms in $H^2 (X, \R )$ corresponding to some rational linear space $\Lambda \subset 
H^2 (X, \R )$.  By choosing $\rho -1$ integral classes $D\in Q$ which together with  $E$ span 
$H^2(X, \R )$, we may assume that 
either the above recipe yields some big movable class, maybe taking a 
sum of two movable classes to obtain bigness, 
or all of them yield movable classes in a (known)  linear space $\Lambda$ of codimension one obtained as above, corresponding on some minimal model to an elliptic fibre space structure.   In the latter case, let $l$ denote a primitive integral form  
defining the hyperplane $\Lambda$; our construction ensures that we may assume that $l>0$ on $Q$ and $l(E) >0$.

In this case, using the results from Section 4 of \cite{WilBd}, we may suppose that $D$ gives rise to $mD = \Delta +rE$ with $\Delta \in \Lambda$ movable, 
and that 
$(m+1) D = \Delta _1 + r_1 E$ with $\Delta_1 \in \Lambda$ movable.  Thus $m l(D) = r l(E)$ and $(m+1) l(D) = r_1 l(E)$.  
Since $(m, m+1) =1$, this implies that $l(E)$ divides $l(D)$.  Now choose an arbitrary integral class $F$ with 
$l(F) = 1$ and (with $D$ as above) a positive integer $N$ with $D' = ND + F \in Q$.  Since $l(E)$ does 
not divide $l(D')$, the above recipe now yields a movable class $\Delta ' \not\in \Lambda$.  With the movable class 
$\Delta$ defining a map to a surface, the sum $\Delta +\Delta '$ is movable and big.

Summing up therefore, given knowledge of the linear and cubic forms on $H^2 (X, \Z )$ and the class of 
the (unique) rigid non-movable surface $E$, 
we can find a finite set of integral classes, at least 
one of which must be movable and big.
\end{proof}\end{prop}

\begin{rem} When the Picard number is odd and the Hessian at $E$ is strictly positive, there is a short-cut in the first part of this Proposition.  In this case then the Hessian at $-E$ it is strictly negative; 
since the Hessian at $D - \lambda E$ is non-negative, and so we obtain a bound $\lambda _0$ as claimed.  \end{rem}
Having proved Proposition 3.3, 
we can then apply Proposition 1.1 from \cite{WilBd} to deduce boundedness once we know the linear and cubic forms on $H^2 (X, \Z)$ and the class $E$.

If $c_2 \cdot E \le 0$, then we have only a bounded number of possibilities for the values $c_2\cdot E$ and $E^3$ by Proposition 2.2 from \cite{WilBd}.  We now prove the same fact 
when $c_2\cdot E > 0$.  From Corollary 3.2 we have $c_2 \ge 0$ on $P$, and the inequality fails to be strict (away from the origin) only if the hyperplane $c_2 =0$ is tangent to $P$ along some ray also lying in $\bar Q$.  

We assume therefore that $c_2 \cdot E > 0$, and the linear form $c_2$  is non-negative on 
$P$, 
where $P^\circ$ is the cone on the bounded component of the cubic, which we have shown above must 
contain the \Kahler cone.  By taking double polars (viz. $D \mapsto D^2$) we get a differentiable map $P \to P^*$, the dual cone in 
the dual space, for which the boundary (a cone on an $S^{\rho -2}$) maps surjectively to the boundary of $P^*$.  Thus the map $P \to P^*$ is also surjective.  
Moreover, this map is injective 
on $P$, since if $D_1 ^2 \equiv D_2^2$, then $(D_1 - D_2)(D_1 + D_2) \equiv 0$.  But $P$ is both convex and contained in the index cone, and so 
in particular it follows that the Hessian is non-vanishing at $D_1 + D_2$.  Therefore the map $L \mapsto (D_1 + D_2)\cdot L$ is injective from the real vector 
space to its dual, and so $D_1 = D_2$.  Thus for any linear form $l$ (for instance $l = c_2$) which is non-negative on $P$, we have  
a unique real class $D\in P$ for which $l = D^2$; if $l$ is strictly positive on the non-zero classes in $P$, then $D\in P^\circ$.  
\begin{rem}
For a smooth ternary cubic, we consider the corresponding elliptic curve in $\P^2 (\R )$.
There will be four (some possibly complex, or maybe coincident) poles 
corresponding to a projective  line, namely in our case $c_2 =0$; these are obtained as the four base-points of the pencil 
of polar conics corresponding to the points of the line, and there is a similar result for higher $\rho$ where we are 
looking for the base points of a higher dimensional linear system of quadrics.
  Of these an even number will be real and the complex poles occur as complex conjugates; some of the real poles may however yield $D^2 = -c_2$ rather than $c_2$  and so will have $D^2$ negative on $P^\circ$.
In  the case $\rho =3$, and the line $c_2 =0$ intersects the Hessian curve  in three real points, each such point $A$ gives rise to a reducible conic in $\P^2 (\R )$ (a line pair), one of the lines of which is 
 the line joining the two projective  points of the visible extremity of $P$ from $A$.  
 Each of the real intersection points of the line with the Hessian curve gives rise to a projective line through 
 the unique pole $D\in P^\circ$ with $D^2 =c_2$, and so $D$ is explicitly determined by the geometry.  
 \end{rem}
 
  \begin{prop} Under the assumption that the form $c_2$ is non-negative on $P$, the closed positive cone on the bounded component of the cubic hypersurface, 
  knowing the cubic and linear forms for $X$ 
  will also bound $c_2\cdot E$ and $E^3$ for the unique rigid non-movable class $E$.
  \begin{proof}We may assume that $c_2\cdot E>0$; we commented above that either $c_2$ is strictly positive on the non-zero classes in $P$, or the hyperplane $c_2 =0$ is tangent to $P$ along some ray generated by a class in the boundary.  
  
  We saw above that when $c_2$ is strictly positive on the non-zero points of $P$, there exists a unique real point $D \in P^\circ$ with
  $c_2 = D^2$.  In the special case where $c_2 =0$ is tangent along some ray in the boundary of $P$, we have that $c_2 = D^2$ for a unique $D$ in that ray.
    
  Choose linearly independent rational points $A_1 , A_2$
  on the projective hyperplane $c_2 =0$, whose polars are rational projective quadrics $Q_i =0$ ($i=1, 2$), both containing the unique point $D \in P$ with
  $c_2 = D^2$.  Choose a rational point $B$, and choose rational values $a_1 , a_2$ such that the quadric $S = a_1Q_1 + a_2 Q_2$ vanishes at $B$.  Letting $A =  a_1A_1 + a_2 A_2$, we obtain a rational point on the projective  hyperplane $c_2 =0$, giving rise to a rational projective quadric $S$ through 
  $D$, whose corresponding affine cone will be denoted by $S^*$.  By choosing $B$ appropriately, we may assume that $A$ does not lie on the Hessian (and in the special case, is not in the ray generated by $D$) and then 
  the projective quadric $S$ is smooth, and so in particular its set of rational points is dense.  
  The space of homogeneous linear forms which vanish at $A$ is 
   a $\rho - 1$ dimensional space $V$ containing $f_1 = c_2$, and we can extend this to a basis 
   $f_1 , \ldots f_{\rho -1}$ of $V$ with the $f_i$  ($i>1$) close to $f_1$ and strictly positive on the non-zero classes in $P$.  These 
   additional linear forms then correspond to 
   points $D_2 , \ldots , D_{\rho -1}$  in $P^\circ$ near $D$ on $S=0$ for which $D_i^2$ is a multiple of $f_i$ for all $i >0$.   
   We may perturb $ f_2 , \ldots , f_{\rho -1} $ in $V$ (retaining strict positivity on $P$) and we are also allowed to scale the $f_i$; using the density of rational points on $S$,
   in this way we may assume that we still have a basis but with $f_i$ rational multiples of $D_i ^2$ for all $i>1$ with 
   the $D_i$ \it integral \rm points in $S^*\cap P^\circ$, which may be taken to be primitive, whose corresponding points in the projective plane are near the point given by $D$.  The 
   intersection of the 
   all the projective  hyperplanes $f_i =0$ in $\P^{\rho -1}$ for all $i\ge 1$  is just the point $A$.  
   For each $i >1$, we let $g_i =0$ define a hyperplane in the pencil of hyperplanes containing 
   $f_1 =0$ (recall $f_1 = c_2$) and $f_i =0$, with $g_i$ also close to $f_1$, 
    such that the open set $\{ f_i >0 \} \cup \{ g_i >0 \}$ 
   contains all of the closed half-space $c_2 \ge 0$ apart from the base locus of the pencil, 
 namely the codimension two linear space given by the vanishing of $c_2$ and $f_i$.  Thus the union of the open sets $\{ f_i >0 \} \cup \{ g_i >0 \}$  for $i>1$ covers 
 all of the closed half-space $c_2 \ge 0$ apart from $A$.
 
 We now parallel translate the hyperplanes $g_i =0$ to another point $A'$ of the hyperplane $c_2 =0$ where $g_i (A') <0$ for all $i>1$, so that the affine linear forms for the 
 translated hyperplanes are $ h_i = g_i + \varepsilon _i$ for suitable $\varepsilon _i >0$. 
 
\begin{claim}We may assume that these translated hyperplanes only intersect $P$ at the origin.
 \end{claim}
 
 In the general case (where $c_2$ is strictly positive on the non-zero points of $P$), this can be achieved by choosing the original $g_i$ sufficiently close to $c_2$ so that none of the $g_i$ vanish on the non-zero points of $P$, and then taking $A' $ in the 
 hyperplane $c_2 =0$ close enough to $A$ so that this property is preserved.  In the special case,  all the hyperplanes $g_i =0$ will intersect $P^\circ$ non-trivially, and so the argument doesn't work.  Instead we need to take an appropriate $A'$ on the affine line through $A$ and $D$ (in particular with $D$ on the line segment $AA'$) with the required property.
 
 In both cases, we suppose that the translated hyperplanes are given by linear forms $h_2, \ldots , h_{\rho -1}$, and then 
 we have a second desired property that the open half-spaces $f_i >0$ and $h_i >0$ for $i= 2, \ldots , \rho -1$ together  cover the closed half-space $c_2 \ge 0$.  Note that both the properties of the hyperplanes being disjoint from $P$ and the $2(\rho -2)$ open half-spaces covering the closed half-plane $c_2 \ge 0$, are \it open \rm properties.  So we may perturb the point $A'$ on $c_2=0$ by a small amount and retain these properties.  So in particular, we may now take $A'$ rational on the hyperplane and such that the corresponding quadric projective hypersurface $S'$ is smooth with rational points being dense.  The $h_i$ can be written uniquely as $G_i ^2$ with the $G_i$ classes on $(S')^* \cap P^\circ$, 
 where $(S')^*$ is the affine cone of the quadric hypersurface $S'$, 
 but by a small perturbation of the $h_i$ we may assume as above, at the expense of the $h_i$ being positive rational multiples of the $G_i^2$,   that the $G_i$ are integral primitive points of $(S')^* \cap P^\circ$.
 
   We may therefore run the above algorithm to find primitive integral points
$D_2, \ldots , D_{\rho -1}, G_2, \ldots G_{\rho -1}$ in $P^\circ$, at least one of which must be in the open 
subcone $Q$, since we are in the case where $c_2 \cdot E >0$.

 Thus we may assume that an integral $D \in Q$ has now been given, and hence for some multiple depending on $D$ we have 
    that $h^0(X, \O_X (mD))>1$, and so $|mD | = |\Delta | + rE$, with some multiple of $\Delta$ 
mobile (and hence $\Delta$ movable), and $r$ a non-negative integer.  Given that $c_2 \cdot E > 0$ 
and $c_2\cdot \Delta >0$, we either have $r=0$  or we
have  obvious upper bounds on both $r$ and 
 $c_2 \cdot E$ in terms of $D$, the latter 
also giving a lower bound on $E^3$ by Proposition 2.2 of \cite{WilBd}. 
If $r=0$, then we have already found a big movable integral divisor $D$
and so boundedness for the Calabi--Yau threefolds follows from Proposition 1.1 from  \cite{WilBd}, and the Proposition follows in this case too.
 \end{proof}\end{prop}

 The plan now is to apply the results of \cite{WilBd} to obtain boundedness.  
 We have shown that $c_2\cdot E$ and $E^3$ only take a finite set of possible values, using Proposition 2.3 of \cite{WilBd} and Proposition 3.6 above.
 If this information were to bound the possible classes of $E$, then Proposition 3.4 would yield a finite set of classes, one of which would be movable and big,  and hence  the results of \cite{WilBd} yield boundedness,
 We therefore need a result which says that, given the cubic and linear forms on $H^2(X, \Z )$, there 
are only finitely many possible integral classes $E$ with given values of $E^3$ and $c_2 \cdot E$.

 \begin{lem}Suppose the Picard number $\rho (X)=3$ and the cubic form defines a real elliptic curve with two components, where $X$ contains only one rigid non-movable surface;
 then for any given values of $E^3$ and $c_2 \cdot E$, there are only finitely many possible classes for the rigid non-movable surface $E$.
 
 \begin{proof} As usual, we assume that $X$ is general in moduli. Crucial to this result are the real points of intersection of the  projective line $c_2 =0$ with the real elliptic curve.
 We let $S = S^2 \subset \R^3$ be the 
 unit sphere given by $x^2 + y^2 +z^2 =1$, so for each point of intersection there are two (antipodal)  points of $S$; unless the point of intersection is an inflexion 
 point, we need only consider the corresponding point $A\in S$ for which $H(A) >0$.  In this case, we take open  neighbourhoods $A \in N_1 \subset N_2 \subset S$ of 
 $A$, with say $N_i$ being a spherical disc of radius $\epsilon _i$, where $ \epsilon _2 = 2\epsilon _1$.  We assume that $\epsilon _1>0$ has been chosen sufficiently 
 small so that $H>0$ on $N_2$ and the convex subset of $P^\circ$ given by $E\cdot D^2 > 0$ for all $E \in N_1$ is non-empty.  We choose a fixed integral class $D$ in this open cone.  Thus, if the class of a rigid non-movable surface $E$ lies in the cone on $N_1$, then $h^0 (X, \O_X (mD)) > 1$ for some fixed $m>0$ depending on our choice of $D$ but 
 not on the class $E$. 
 I claim that there can only be finitely many 
 possible classes $E$ of rigid non-movable surfaces in the open cone on $N_1$.  This follows since there exists $M>0$ such that if $E$ is in the closed cone on $N_1$ with $||E|| >M$, then $mD -r E \in -N_2$ for all $r >0$, giving a contradiction since $H<0$ on all such classes (and so none could be movable).  Thus for any possible class of 
 a rigid non-movable surface $E$ in the cone on $N_1$, we have $||E|| \le M$, and hence there are only finitely many such classes.
 
 We now consider the case where the line $c_2 =0$ intersects the cubic at an inflexion point, and we consider the corresponding points $\pm A \in S$, where without loss of generality we may take $A = (0,1,0)$.  So the conic defined by $A\cdot D^2 =0$ consists of two lines intersecting at $Q_1$ on the Hessian, one of them not intersecting 
 the projectivised cone $P$ and the other joining the two points on the boundary of the projectivised cone where the tangent lines are vertical (the points of the visible extremity from 
 $A$).  Moreover direct calculation verifies that the linear form $E^2$ is negative on that part of the affine plane $z=1$ given by $x> -1/(k-1)$, and in particular therefore $E^2$ is 
 strictly negative on the non-zero elements of $P$.  We can now find an open neighbourhood $A\in N_1 \subset S$ with closure $\bar N_1$ such that 
 $E^2$ is strictly negative on the non-zero elements of $P$ for all $E \in \bar N_1 \subset S$ and such that the convex subset of $P^\circ$ given by $E\cdot D^2 > 0$ for all $E \in N_1$ is non-empty.  We choose a fixed integral class $D$ in this open cone.  Thus, if the class of a rigid non-movable surface $E$ lies in the cone on $N_1$, then $h^0 (X, \O_X (mD)) > 1$ for some fixed $m>0$ depending on our choice of $D$ but 
 not on the class $E$. I claim that there can only be finitely many 
 possible classes $E$ of rigid non-movable surfaces in the open cone on $N_1$.  This follows since there exists $M>0$ such that if $E$ is in the closed cone on $\bar N_1$ with $||E|| >M$, then $(mD -rE)^2 = m^2D^2 -2rmD\cdot E + r^2E^2 < 0$ on the non-zero elements of $P$ for all $r>0$.  Since then $(mD-rE)^2\cdot L < 0$ for any ample class $L$, we deduce that $mD-rE$ cannot be movable for any $r>0$.  Thus for any possible class of 
 a rigid non-movable surface $E$ in the cone on $N_1$, we have $||E|| \le M$, and hence there are only finitely many such classes.
 
 Suppose now there were infinitely many possible classes $E$ with given values of $E^3$ and $c_2\cdot E$, 
 for the rigid non-movable surface, then $||E||$ would be unbounded.  We could then find a sequence of distinct classes $E_i$, with corresponding classes $ \bar E_i := E_i/||E_i|| \in S$ tending to one of the points in $S$ corresponding to an intersection point of the line $c_2 =0$ and the cubic.  Therefore, for one of the neighbourhoods $N_1 \subset S$ we found above, we have $\bar E_i \in N_1$ for all $i\gg 0$, contradicting the previous conclusions.
\end{proof}\end{lem}

   For $\rho >3$ we  do not at present have available any similar general result to Lemma 3.7, and in fact one should only hope for a finiteness result modulo isomorphisms 
 of $H^2 (X, \Z )$ fixing the linear and cubic forms, which would suffice for our purposes.  
   We have however proved all the previous results in this section for the general case, and in particular for all $\rho >2$, 
 rather giving the simpler proofs for $\rho =3$, to emphasise that it is precisely here and nowhere else that the proofs 
 for the general case fail, in particular of boundedness when $\rho >3$ and there are two real components of the cubic.  
  When dealing with the cubic having one real component, as studied when $\rho =3$ in the next section, we will also be confronted with the fact that the geometry of real cubic hypersurfaces is just far more complicated if $\rho >3$.
 
In summary however, for $\rho =3$ we have proved the following statement.

\begin{cor} Suppose that $X$ is a \CY threefold with $\rho =3$, with given  linear and cubic forms on $H^2 (X, \Z )$, where the cubic defines a real  elliptic curve with two components.  
  If $X$ contains at most one rigid non-movable surface, then 
$X$ lies in a bounded family.
\end{cor}

\section{Case when the elliptic curve has one component}

  We now wish  to study the case when the elliptic curve has only one real component, and so if the Hessian is smooth, 
it has two components.  Throughout the section, we may without loss of generality 
assume that $X$ is general in moduli.  The four components of the positive index cone are  as 
described in Section 1.  Recall that there are two special values for $k<1$ where changes occur, namely $k =0$ an $-2$.  Away from these two values (which we shall consider later in this Section), we wish to describe the Steinian map $\alpha$.

  If as before the inflexion points of the 
cubic (and hence of the Hessian) are denoted $B_1, B_2, B_3$, then the tangents there  to 
the cubic $F$ (which we saw are just the asymptotes to the three affine branches) will be 
tangent to the Hessian at three distinct points $Q_1, Q_2, Q_3$.  Having chosen one of the $B_i$ as the zero of the group law, the corresponding point $Q_i$  where the tangent to $F$ at $B_i$ is tangent to the Hessian 
is just one of the 2-torsion points of the Hessian.   Moreover the second polar of $F$ with respect to each of these three points $B_i$ will be the tangent to the Hessian at the corresponding point $Q_i$.  

It is easiest to understand what is going on dynamically.  
For $k>1$, we found a precise description of $\alpha$, 
where the tangent to $F$ at the each $B_i$ is tangent to an affine branch of the unique connected component of $H$.  As $k \to 1$, the affine branches of both 
the real curves given by $F$ and $H$ tend to segments of the line at infinity between the relevant inflexion points, 
whilst the bounded component shrinks to a point, so that for $k=1$ both $F$ and $H$ define the line at infinity plus an isolated point at the centroid $({1\over 3}: {1\over 3}:1)$ of the triangle of reference.  Deforming away from $k=1$ towards zero, the singular  point then 
expands to become the 
bounded component of the Hessian, and the line segments at infinity that were limits as $k\to 1+$ 
of the affine branches of $F$ deform to affine branches of $H$ and  the line segments that were limits as $k\to 1+$ 
of the affine branches of $H$ deform to affine branches of $F$.  Recall here that $H_k = -54k^2 F_{k'}$ from Section 1, and so one does expect the regions occupied 
by the unbounded affine  branches of $F$ and $H$ to switch over.
 In particular, the tangents to the cubic at each inflexion point are tangents to the unbounded affine branches of $H$, similar therefore to the case $k>1$.  
Thus in this case, the Steinian map $\alpha$ sends each inflexion point $B_i$ to a point on the unbounded component of $H$, and hence gives an involution on both 
connected components of $H=0$.  

The next change occurs at $k=0$, where the bounded component together with the three affine branches of $H$ just tend to the three
real lines determined by the triangle of reference.  To see what happens to the Steinian map, it is probably easiest to look 
at the value $k= -2$; here the Hessian is just the line at infinity together with the isolated point $({1\over 3}: 
{1\over 3} :1)$ and all three asymptotes pass through this point.  As one deforms in either direction away from $k=-2$, this point expands to give the bounded component of the Hessian and each asymptote of $F$ will 
now be tangent to the bounded component of the Hessian, which by continuity will also be the case for all 
$-2 < k <0$.  Thus for $k < -2$ and $-2 < k < 0$, the Steinian map $\alpha$ interchanges the two components of $H=0$.  The bounded component will be contained in (and tangent to) the asymptotic triangle given by the lines $x ={1\over {1-k}}$, $y = {1\over {1-k}}$ and $x + y = {{k}\over {k-1}}$.  For $k >-2$, the asymptotic  triangle will be given by the inequalities $x \le{1\over {1-k}}$, $y \le {1\over {1-k}}$ and $x + y \ge {{k}\over {k-1}}$, whilst for $k<-2$ it will be given by $x \ge {1\over {1-k}}$, $y \ge {1\over {1-k}}$ and $x + y \le {{k}\over {k-1}}$.

What is occurring here is that for each value of $k' > 1$, there are three possible values of $k$ for which $H_k$ is a multiple of $F_{k'}$, one with $k<-2$, one with $-2 < k < 0$ and one with $0< k <1$.  If we choose an inflexion point  $B_3$ say, the tangent to $F_k$ at $B_3$ will be tangent at 
one of the three 2-torsion points of $F_{k'} = H_k$, the one on the unbounded component if $0< k <1$, and the ones on the bounded component in the other two cases.   

For any given point $A$ in the affine plane $z=1$, we let $G_A$ denote the function of $A = (a,b,1)$ given by $A\cdot D^2$, and we wish to understand how 
$G_A =0$ intersects not only the affine branches of $F$ but also the unbounded affine branches of the Hessian.   We will therefore need to understand this in all the three cases detailed above, as the Steinian map will be different in the three cases.

\begin{prop} Suppose that $X$ is a  \CY threefold with Picard number $\rho =3$ whose corresponding real cubic curve $F_k$ is smooth  with one real component,  with the Hessian also smooth, and that there is at most one   
rigid non-movable surface $E$ on $X$.  Then such a surface exists and its class $E$ lies in the closure of the half-cone on which $H\ge 0$ determined by 
 the bounded component of the Hessian, and moreover $X$  lies in a bounded family.

\begin{proof} We recall (Proposition 1.1) that the component $P^\circ $ of the positive index cone containing the \Kahler cone would have to be hybrid.  We now let 
 $C_1$ denote the unbounded branch of $F=0$ which lies in the region 
$x>0, y>0$ and $x+y >1$, and $C_2$ the unbounded branch of the Hessian lying in the sector $x<0, y<0$.  
There is then a hybrid component $P^\circ$ of the positive index cone whose boundary consists of the positive cone on $C_1$ together with the negative cone on $C_2$, the two parts meeting along rays generated by $(0,1,0)$ and $(1,0,0)$, and without loss of generality we may assume that this is the component which contains the \Kahler cone.

As observed in the Introduction, there will be a (unique) rigid non-movable 
surface class $E$ on $X$.
In the special case where $E$ lies in the closure of the half-cone on the bounded component of the Hessian on which $H\ge 0$, 
we do however have $E^3 >0$ and hence as observed in  \cite{WilBd} that $1\le E^3 \le 9$.  
 The fact that the cubic is strictly positive on the non-zero elements of this closed half-cone, 
in this case implies that there will 
only be finitely many classes with the given possible values of $E^3$.
For each such class $E$, the condition $G_E >0$ will define an open subcone of $P^\circ$ (in fact by 
considering the zeros of $G_E$ on the projectivised boundary of $P$, arguing as below 
the subcone is in all cases connected, and unless the cubic is $F_k$ with $k<-2$ will be all of $P^\circ$); without however using 
these facts, for any component $Q$  of this open cone
we may choose an integral $D \in Q$; 
if $Q$ also contains the \Kahler cone,   then by Proposition 4.1 and Lemma 4.3 of \cite{WilBd} we can find 
 a positive integer $m$ depending on $D$ for which $h^0(X, \O_X (mD)) >1$.
 For $A = -E$,  it is however clear that $D +\mu A$ is in the other half-cone 
on the bounded component of the Hessian for $\mu \gg 0$ and hence has value of the Hessian negative.  Thus there is an upper bound $\lambda _0$ for $\lambda >0$ for which $D - \lambda E$ can be moveable, and so we find a finite set of classes, at least one of which must have a multiple which is mobile.  
The class so found might not be big, but by running the same argument with a different choice of $D' \in Q$ not in the plane spanned by $D$ and $E$, we can as argued in the proof of Proposition 3.3  ensure that it is either big or corresponds to an elliptic fibration on an appropriate minimal model for $X$;  as further argued in that proof we may then find a finite set of integral classes with at least one being movable and big, 
and thus by Proposition 1.1 from \cite{WilBd} 
the claimed result follows. 

From now on therefore, we shall assume also that 
$E$ does not lie in the closure of the half-cone on the bounded component of the Hessian on which $H\ge 0$.
We are interested in integral classes $E$ with $H(E) \ge 0$; we let $A$ denote the class
 $(a,b,1)$  in the affine plane $z=1$ given by its intersection with the line through the origin and $E$ --- 
 we shall comment later that the proofs extend to the case when $E$ lies in the plane
 $z=0$.

Using Lemma 3.3 of \cite{WilBd}, a connected component of the subcone given by $A\cdot D^2 > 0$ 
 will  be convex  if $H(A) \ge 0$, for instance when $A$ is a positive multiple of $E$,
 whilst a connected component of the subcone given by $A\cdot D^2 < 0$ 
 will  be convex  if $H(A) \le 0$, for instance when $A$ is a negative multiple of $E$.
 Let $Q$ denote a connected component of 
the subcone of $P^\circ$ given by $E\cdot D^2 >0$, therefore convex.  In the cases when $0<k<1$ or 
$-2 < k < 0$, we shall see that the conic $G_A =0$ cuts 
the projectivised boundary of $P$ in at most two points and so the subcone of $P^\circ$ given by 
$E\cdot D^2 >0$ is itself connected; thus for a given class $E$, the cone $Q$ is unique;
we shall see however that this will no longer always be true when $k<-2$.

\begin{claim} Under the above assumptions, if $A$ is a positive multiple of $E$, part of the projectivised boundary of $Q$ corresponds to an open arc of $C_2$  consisting of points that are visible  from $A$. 
If $A$ is a negative multiple of $E$, part of the projectivised boundary of $Q$ corresponds to an arc of $C_2$ that is not visible  from $A$.
\end{claim}
Given this claim, since  the corresponding  part of the boundary of $P$  consists of rays on the negative of points in the given arc on $C_2$,  
the part of the boundary of $Q$  in question is not visible from $E$.  Assuming that $E$ does not lie on the line $z=0$ at infinity, 
this shows via Proposition 1.2 that the \Kahler cone cannot be contained 
in a hybrid component of the positive index cone, and in the light of 
Proposition 1.2 the required result then follows.

For the case when the cubic is $F_k$ with $0<k<1$, the argument is essentially identical to the argument for $k>1$ given in Section 2, modulo the fact that the 
regions of the plane occupied by the unbounded branches of $F$ and $H$ have switched over.    For $k>1$, 
the bounded component of $F$ essentially plays no role in the argument, and for $0 < k <1$, 
the bounded component of $H$ essentially plays no role in the proof of the above Claim.  The case when $E$ lies on the line at infinity also follows as in the argument from Section 2.

The other two cases (where the Steinian map interchanges the two components of $H$) will be very similar to each other;  we shall first give the argument 
when $-2 < k < 0$ to show that under the above assumptions the \Kahler cone 
of the \CY threefold $X$ cannot 
lie in a hybrid component of the positive index cone, 
 and afterwards we shall comment what is different in the case $k<-2$.  

We note in the first of the two cases that the tangent line to $F$ at $B_3$ is the line 
$x+y = e_2 = {k\over {k-1}} = {{-k}\over {1-k}}$, and this is tangent to  the bounded component of the Hessian at the point $(e_2/2 , e_2/2)$.  
If we take $B_3$ as the zero of the group law, this is just a 2-torsion point of $H$.  Moreover the bounded component of the Hessian
 is inscribed in the asymptotic triangle, and so in this case it lies above (and touches) the line $x+y = {k\over {k-1}}$.  
Also playing a role will be the other other two lines through $B_3$ that are tangent to the Hessian and yielding 2-torsion points of the Hessian; these have the form $x+y = e_3$, 
corresponding to the other tangent to the bounded component  and $x+y = e_1$ corresponding to the tangent to the unbounded component of $H$.  
Explicitly, if the Hessian is  (up to a multiple) the Hessian of $F_{k_i}$, 
where $0 < k_1 <1$, $-2 < k_2 = k  <0$ and $k_3 < -2$, then $e_i = k_i/(k_i -1)$; moreover $e_1  < e_2 < e_3$.

We have the Steinian involution on the Hessian which in these two cases interchanges the two components; explicitly we 
let $Q_1$ be the point on the bounded component of the Hessian whose tangent is also the tangent to $F$ 
at $B_1$, and $Q_2, Q_3$ defined similarly with respect to the inflexion points $B_2, B_3$.  The latter we saw 
above was the point $(e_2/2 , e_2/2)$ and the tangent line $x+y =e_2$, where $e_2 = {k\over {k-1}}$.

Under the Steinian map, the branch $C_2$ (going from $B_1$ to $B_2$)  
of the Hessian corresponds to the arc (not containing $Q_3$) on the bounded component going from $Q_1$ to $Q_2$.  We now argue similarly to Section 1.  For a given class $A = (a,b,1)$ not in $P$,  it is clear how many times the quadric $G_A = A\cdot D^2$ cuts (the closure of) $C_1$ --- it will cut it twice if $a\ge {1\over {1-k}}$ and $b\ge {1\over {1-k}}$ (since there will be two tangents to $C_1$),  it will not cut $C_1$ at 
all if $a< {1\over {1-k}}$ and $b< {1\over {1-k}}$, and will cut it once in the other cases.  We now ask how many times and where the quadric cuts $C_2$.  To answer this question, we are looking for the tangents from $A$ to 
the arc (not containing $Q_3$)  from $Q_1$ to $Q_2$ 
on the bounded component of the Hessian.  Here the answer is twice (with multiplicity) if 
$A$ is in the region bounded by $a 
= {1\over {1-k}}$, $b ={1\over {1-k}}$ and by the arc specified $Q_1 Q_2$, none for any  other points with 
$a< {1\over {1-k}}$ and $b< {1\over {1-k}}$ or with $a> {1\over {1-k}}$ and $b> {1\over {1-k}}$, and 
precisely once otherwise as there is exactly one tangent to the given arc $Q_1 Q_2$.  
Moreover, the midpoint of 
the arc is the point $(e_3/2 , e_3/2)$ noted before, and under the Steinian map this point corresponds to the midpoint 
of $C_2$, namely the intersection of $C_2$ with $x=y$.  So being more precise still, $G_A$  will cut $C_2$ in the part given by $y\le x$ if and only if $a \le {1\over {1-k}}$, $b \ge {1\over {1-k}}$ and $a+b \ge 
e_3$ \it or \rm $a \ge {1\over {1-k}}$, $b \le {1\over {1-k}}$ and $a+b \le 
e_3$.

The Claim is true for all the possibilities detailed above when $-2<k<0$ for where the function $G_A$ vanishes on $C_1 \cup C_2$.
 We first explain the proof of the claim in detail 
  in the case 
$a\le {1\over {1-k}}$, $b\ge {1\over {1-k}}$.

We assume first that $A$ is a positive multiple of $E$.
Of course if say we have $b= {1\over {1-k}}$ with  $a$ arbitrary, we have that the line $y = {1\over {1-k}}$ is both tangent to $C_1$ at $B_2$ and to the bounded component of the Hessian at $Q_2$, and so $G_A$ has 
just a double zero at $B_2$.  It follows that either $G_A  <0$ or $G_A >0$ on all of $P^\circ$, 
yielding either that $Q$ is empty  or $Q= P^\circ$ --- in the latter case clearly there exists an arc of $C_2$ 
which is visible from $A$.
We assume therefore that we have strict inequalities for $a$ and $b$ above.  There is precisely one point $U$ on the open arc $Q_1 Q_2$ of the 
specified bounded component of the Hessian for which $A$ lies on the tangent line, an so the 
argument via the Steinian involution shows that there is precisely one point, namely $\alpha (U)$ on $C_2$ at which $G_A$ vanishes.  A mechanical calculation verifies 
that $G_A (B_1) < 2k/(1-k) <0$, and hence we must have 
$G_A(B_2 ) >0$.  Thus there is an arc in $C_2$ (namely near $B_2$) which is part of the projectivised boundary of $Q$ and which is visible from $A$.

Assume now  that $A$ is a negative multiple of $E$, and hence $H(A) \le 0$.  
We know that $E\not \in P$ and so we have a stronger inequality that $a < -{1\over {k'-1}}$, where $x = -{1\over {k'-1}}$ is an asymptote for $C_2$, namely the tangent to the Hessian at the inflexion point $B_1$.  
We will still deduce as above that $G_A (B_1) <0$, and thus there is an arc of points of $C_2$ on 
which $G_A$ is strictly negative and which are not visible from $A$ (namely near $B_1$), noting that $B_1$ itself is not visible from $A$
because of the stronger  inequality on $a$.  

The second case we explicitly check is when $A$ is in the open region bounded by the arc (not containing $Q_3$) of the bounded component of the Hessian from $Q_1$ to $Q_2$, and  line segments from 
the lines $x={1\over {1-k}}$ and $y={1\over {1-k}}$ --- in this case the arc lies in the quadrant $x \le {1\over {1-k}},\  y \le {1\over {1-k}}$.  Note that $H(A) >0$ and so $A$ is a positive multiple of $E$.  One checks that the negative of 
 the arc of $C_2$ from $\alpha (Q_1)$ to $\alpha (Q_2)$ corresponds to part of the boundary of $Q$, and this arc is clearly visible from $A$.  If now we take $A$ on the 
 arc $Q_1 Q_2$ (not containing $Q_3$) of the bounded component of the Hessian, then $G_A$ is a line pair only intersecting $C_2$ at one point and not intersecting 
 the affine branch $C_1$ at all, and we check that $G_A =0$ is a pair of complex lines whose only 
 real point is the  intersection point with $C_2$ --- recall that the end points of the arc correspond to the inflexion points $B_1$ and $B_2$.  Thus according to whether $A$ is a positive or negative multiple of $E$, we have that $G_A$ is negative on $P^\circ$ or positive on $P^\circ$.  The former case clearly does not occur and in the latter case all of $C_2$ corresponds to boundary points of $P$.

The reader should check that the Claim remains true for the other possibilities for $A$ when $-2 < k <0$; these are either easier or follow by symmetry from the first case considered.

Let us now consider  the above claim in the remaining case with smooth Hessian, namely $k<-2$; we will see that the claim continues to hold.  As before, we may assume that 
$E$ does not lie in the closure of the half-cone on the bounded component of the Hessian on which $H\ge 0$,
and does not lie in the plane $z=0$;
we let $A=(a,b,1)$ denote the point in the affine plane $z=1$ determined by $E$.
  Let $Q$ denote a connected component of the subcone of $P^\circ$ given by $E\cdot D^2 >0$, 
by Lemma 3.3 of \cite{WilBd}, a convex subcone of $P^\circ$

As before it is clear for $A = (a,b,1)$ how many times (and where) $G_A$ intersects $C_1$.  It will intersect $C_1$ twice if $a \ge {1\over {1-k}}$ 
and $b \ge {1\over {1-k}}$, at no points if $a < {1\over {1-k}}$ 
and $b < {1\over {1-k}}$, and once otherwise.  For 
$a < {1\over {1-k}}$, $b < {1\over {1-k}}$, there are no zeros of $G_A$ on  $C_2$ and $G_A$ would be negative on all of $P$.  If say 
$a \le {1\over {1-k}}$, $b = {1\over {1-k}}$, then $G_A$ would be negative on $P^\circ$.  If $a > {1\over {1-k}}$, $b = {1\over {1-k}}$ then $G_A$ has zeros on the affine branch $C_1$ and the point $B_2$ at infinity, and the whole of $C_2$ is visible from $A$.

An additional feature compared with the previous case is that $G_A$ can intersect the projectivised boundary of $P$ at two points on $C_1$ and two on $C_2$, and that will happen when $A$ is in the open region with boundary consisting of a segment of the line $x = {1\over {1-k}}$, a segment of the line  
$y = {1\over {1-k}}$, and the arc (not containing $Q_3$) of the bounded component of the Hessian between $Q_1$ and $Q_2$; as before $Q_i$ denotes the point on the bounded component of the Hessian where the asymptote to the cubic at $B_i$ is tangent. In contrast to the previous case, this time the arc 
of the Hessian between $Q_1$ and $Q_2$ lies in the quadrant $x \ge {1\over {1-k}}$, $y\ge {1\over {1-k}}$.
We illustrate this situation with the Figures 4 and 5 when $k=-3$, and $A= (0.28, 0.28, 1)$, which lies in the open region bounded by $x= {1\over 4}$, $y= {1\over 4}$ and the 
arc of the Hessian between $Q_1$ and $Q_2$.

 \begin{figure}
      \includegraphics[width=9cm]{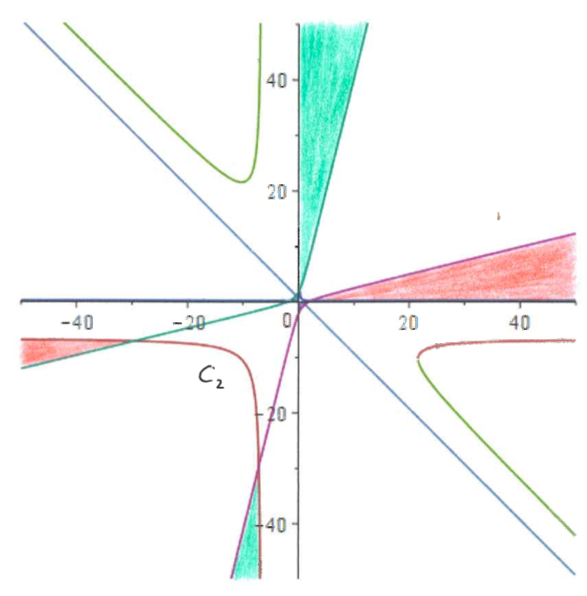}
\caption{Case of four intersections of $G_A$ with projectivised boundary}
\end{figure}   

 \begin{figure}
  \includegraphics[width=9cm]{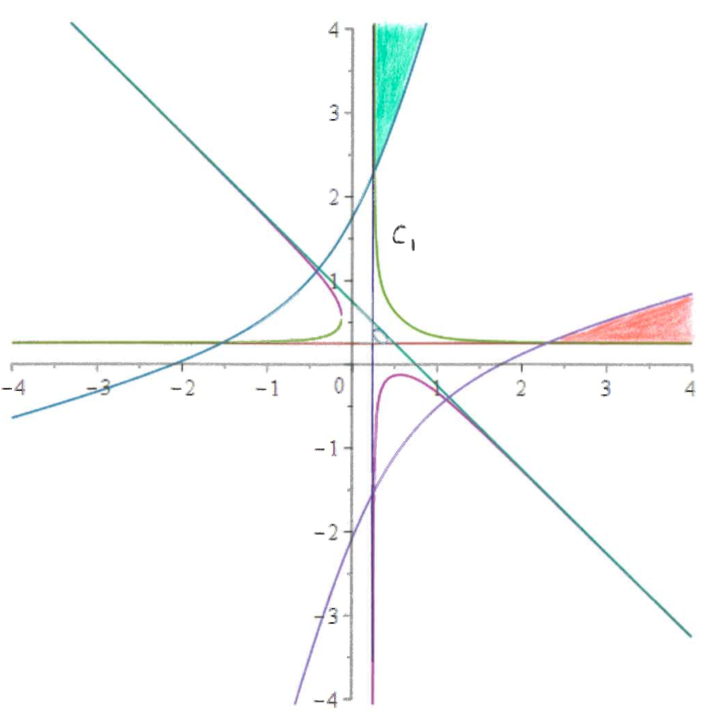}
\caption{Close-up near origin of Figure 4}
\end{figure}

         Figure 5 is just the detail near the origin of the Figure 4, and includes the cubic and its asymptotes, the quadric $G_A$, 
but only the bounded component of the Hessian is visible (only just!), inscribed in the triangle formed by the asymptotes.
      In all such cases $H(A) >0$ and the condition $A\cdot D^2 >0$ defines \it  two \rm open convex subcones of $P^\circ$, but both of these satisfy the Claim made before, 
      that part of the projectivised boundary of the subcone corresponds to an open arc of $C_2$  consisting of points that are visible  from $A$; recall that in this case, the 
      relevant branch $C_2$ of the Hessian is contained in the negative quadrant. 
This is illustrated in the  two pictures, where the \it two \rm open convex subcones of $P^\circ$ given by $A\cdot D^2 >0$ are defined by the regions shaded in red and green.

To justify that these are the correct pictures, let us consider the degenerate case where $A= ({1\over{1-k}} , {1\over{1-k}}, 1)$;  here the conic $G_A =0$  only intersects the projectivised boundary of $P$ at the points $B_1$ and $B_2$.  Moreover, calculation gives the formula (for this choice of $A$) that $G_A (x,y,0) = 3(k+2)xy$, which is negative on the open cone generated by $(1,0,0)$ and $(0,1,0)$.  The index at $A$ is $(1,2)$, and so the open cone given by $G_A >0$ has two (convex) components $U_1$ and $U_2 = -U_1$.  From the above calculation $P \not \subset \bar U_i$  for either $i$, and then we have say $(0,1,0) \in \bar U_1$ and $(1,0,0)\in \bar U_2$.  If now $U_1$ contains a point of $P^\circ$, then 
by convexity of $P^\circ$ and the fact that the points on the open cone generated by $(0,1,0)$ and $(1,0,0)$ are not in $U_1$, there exists a ray of $ \partial U_1$ in $P^\circ$; this 
would imply the existence of another ray in $\partial U_1 \cap \partial  P$, contradicting the previous deduction.
We deduce that $\bar U_1 \cap P$ is the ray generated by $(0,1,0)$ and $\bar U_2 \cap P$ is the ray generated by $(1,0,0)$.  If we perturb $A$ to a point in the interior of the given region, the conic $G_A = 0$ will intersect the projectivised boundary of $P$ in two points near each inflexion point, one on $C_1$ and one on $C_2$, and $G_A >0$ defines two convex cones in $P^\circ$; this is the only case where $P^\circ$ intersects both the cones in $\R^3$ defined by $G_A >0$.  This picture continues to hold true for all $A$ 
in the interior of the region by a continuity argument.  

If $A$ is in the interior of one of the line segments of the region in 
question, then $G_A >0$ will define a single convex subcone of $P^\circ$, although there will be another zero on the projectivised  boundary of $P$ at an inflexion point.   Since 
$A$ is a positive multiple of $E$ and all of $C_2$ is visible from $A$, the above claim continues to hold in this case.
 
If $A$ is a point on the specified arc $Q_1 Q_2$, then $G_A$ degenerates into two real lines intersecting at the point $\alpha (A) \in C_2$, 
each also intersecting $C_1$; so long as $A$ is a positive multiple of $E$, then $P^\circ \cap \{ D\ :\ E\cdot D^2 >0 \}$ has two (convex) components $Q$, and 
the same argument as before proves the claim.  
 
All the other cases are similar to the corresponding cases for $-2 < k < 0$, and so the claim has been verified for affine points $A = (a,b,1)$ for both the cases where the Steinian involution switches the two components of the Hessian; thus part of the boundary of $P$  consists of rays on the negative of points in an open arc in $C_2$,  which are  not visible from $E$.

It is a mechanical check that in both cases  $-2 < k < 0$ and $k<-2$, the 
analogous statement continues to hold if $E$ represents a point on the line at infinity $z=0$.  In summary therefore, assuming that $E$ is not in the closure of the half-cone on 
which $H\ge 0$ determined by the bounded component of the Hessian,
 we deduce for all $E \not \in P$, that there is 
  part of the boundary of $Q$ which is not visible  from $E$ but on which the Hessian vanishes.  
When $H(E)=0$ with $E$ lying on the cone on the unbounded component of the Hessian, it follows as in Section 2 that $E\cdot D^2 >0$ divides $P^\circ$ into convex cones, each of which satisfies the claim.

We can therefore apply Propositions 1.1 and 1.2 to deduce that in the case where the cubic form of a \CY threefold $X$ defines an elliptic curve with one component, with the Hessian also smooth, if there is at most one 
rigid non-movable surface on $X$, then  $X$ lies in a bounded family.
  \end{proof}\end{prop}

Let us now consider the two cases where the Hessian is not smooth, 
namely $k=0$ and $-2$.  We recall that in these cases we also assume that the line $c_2=0$ does not intersect the real elliptic curve at an inflexion point.  
 Here we substitute explicit calculation for the theory of the Steinian map.  As usual, if $E$ denotes a unique rigid non-movable surface on $X$, we let $Q$ denote the component of the subcone of 
$P^\circ$ given by $E\cdot D^2 >0$ containing the \Kahler cone, and we set $B_1 = (0:1:0)$ and $B_2 = (1:0:0)$.

  Let us start with the case $k=0$, where we shall assume without loss of generality that $C_1$ is the branch of $F=0$ in the positive quadrant, with asymptotes $x=1$, $y=1$, 
and $C_2$ is given by the two line segments $L_1$ : $x=0, y\le 0$ and $L_2$ : $y=0$, $x\le 0$, so that the cones on $C_1$ and $-C_2$ meet along the rays generated by 
$(0,1,0)$ and $(1,0,0)$.
 For an affine point $A = (a,b,1)$, one knows when $G_A$ intersects (the closure of) $C_1$, that is when $A$ lies on a tangent line from a point of $C_1$, and one can explicitly calculate when and where 
$G_A = -ax^2 -by^2 -(1-a-b)(1-x-y)^2$ is zero on $C_2$.  
Given the class of $E \not \in P$, one is left with several cases; we leave the reader to check the case when $E$ lies on the line at infinity  $z=0$, and summarise below the various cases where some multiple of $E$ is of the form $A = (a,b,1)$; the actual calculations are tedious but easy and so will be omitted.  In all these cases, we 
would like to argue as in the smooth case, using Proposition 1.2, to obtain a contradiction.  In all cases except for (1) and (2), the quadratic $G_A$ will have one zero on $C_1$ and one on $C_2$.

(1) If $a\ge 1$, $ b\ge 1$, then $G_A$ vanishes twice on (the closure of) $C_1$ and is non-zero on $C_2$.  Moreover $H(A) >0$ and so $A$ is a positive multiple of $E$, and all of $C_2$ is visible from $A$.

(2) Suppose now that $0\le a \le 1$, $0\le b \le 1$ and $A \ne (1,1,1)$.  When $a+b \le 1$, we have that $H(A) \le 0$ and $G_A$ is negative on $P^\circ$.  If $A$ is a positive multiple of $E$, we have that $Q$ is empty and hence a contradiction.  If $A$ is a negative multiple of $E$, then $E^3 >0$ and we deduce that there are only finitely many possible classes for $E$  as in the proof of Proposition 4.1, and hence 
 $X$ lies in a bounded family.
 When $a+b >1$, we show that $G_A$ vanishes once on $L_1$ and once on $L_2$.  We note that $H(A) >0$ and the arc of $C_2$ between the two zeros of $G_A$ is visible from $A$ but not in the visible extremity.

(3) Suppose now that $0 < a < 1$, $b \ge 1$; then $H(A) >0$ and 
$G_A$ has a zero on $L_1$, $G_A(B_1) <0$, $G_A (B_2) >0$ and 
points on $L_2$ are visible from $A$ but not in the visible extremity. 

(4) Suppose $a=0$ and $b> 1$; then $H(A)=0$ and $G_A$ has a zero on $L_1$, $G_A(B_1) <0$, $G_A (B_2) >0$ and both $(0,-1,0)$ and $(-1,0,0)$ are visible from $A$.  The case when $A$ is a positive multiple of $E$ follows as in (3), by taking points on $L_2$  near $B_2$.  When $A$ is a negative multiple of $E$, we could 
consider instead points on $L_1$ near 
$B_1$, but these are also visible from $E$.  We find ourselves therefore left with cases where $H(E) =0$ and $E^3 >0$.

(5) Suppose now $a<0$, $b\ge 1$; then $G_A (B_1 ) <0$ and $G_A(B_2) >0$.  When $a+b <1$, we have $H(A) >0$ and that $G_A$ has a zero on $L_2$.  Note that points on $L_2$ (near 
$B_2$) are visible from $A$.  When $a+b >1$, we have $H(A) <0$ and $G_A$ has a zero on $L_1$, and 
points on $L_1$ (near $B_1$) are not visible from $A$.  If $a+b =1$, then $H(A) =0$ and  $G_A =0$ defines two lines through $(0:0:1)$, one of which only meets $C_2$ at $(0:0:1)$ and the other divides the interior of $C_2$ in two parts.  The fact that 
points on $C_2$ near $B_2$ are visible from $A$ and points on $C_2$ near $B_1$ are not visible from $A$ covers both the cases when $A$ might be a positive or negative multiple of $E$.

(6)  Finally we suppose that $a< 0$ and $0< b < 1$; here we can check that $H(A) >0$ and that $G_A$ is negative on all of $P^\circ$, and so $Q$ would again be empty.

All the other cases with $A = (a,b,1)$ then follow by symmetry.  In all the cases, one can apply Proposition 1.2 to obtain the required contradiction to there being  just one rigid non-movable surface, \it except \rm in case (2) above and $A$ is a negative multiple of $E$,  
\it or \rm
in case (4) above where $H(E)=0$ and $E^3 >0$.  In this latter case we may assume that  $A = (0,b,1)$ with $b > 1$ and $A$ is  a negative multiple of $E$.  Since however $E^3>0$, we know that there are only a finite set of possibilities for the pair of values $E^3$ and $c_2\cdot E$, using Lemma 2.2 of \cite{WilBd}. 

\begin{prop} Suppose that the Picard number $\rho (X) =3$ and the cubic in appropriate real coordinates is the Fermat cubic, and we assume that the elliptic curve does not intersect the 
line $c_2 = 0$ at an inflexion point.  If $X$ contains no more than one rigid non-movable surface, then 
 $X$ lies in a bounded family.

\begin{proof} As observed in the Introduction,  there will be a unique rigid non-movable 
surface class $E$ on $X$.
From the calculations described above and the use of Proposition 2.1, we must be in one of the two exceptional cases detailed in (2) and (4) above.  In the first exceptional case, we saw that there are only finitely possibilities for the class of $E$ and 
the argument from Proposition 4.1 shows that $X$ lies in a bounded family. 

In the second exceptional case, 
$A = (0,b,1)$ with $b> 1$ and $A$ is a negative multiple of $E$.  Hence as we have just argued, $E^3 >0$ and  can only take finitely many values, and the same is then true 
for $c_2 \cdot E$.  Under the assumption that the line $c_2 =0$ does not contain an inflexion point of the elliptic curve, arguing as in 
 the first part of the proof of Lemma 3.7 then shows that 
 there are again  only finitely many possibilities for $E$.  Suppose now that we are given the class $E$, 
and $A= (0,b,1)$ is a negative multiple of it.  Given that $b>  1$, an easy calculation verifies 
the line pair $E\cdot D^2 =0$ does not meet the line segment $L_2$ and meets the line segment $L_1$ in a single point, and that if $-D = (0,y,0)$ then $E\cdot D^2 >0$ for $y\ll 0$ --- the quadric meets $C_1$ also in one point, 
and so the points of $D \in \bar Q$ at which the Hessian vanishes are precisely of the form 
$-D = (0,y,0)$, for $y\ll 0$.  A further elementary calculation verifies 
that $A^2\cdot D'  = -b^2 y - (1-b)^2 (1-y) >0$ for all points $D' = (0,y,1) \in L_1$ with $y\le y_0$ for some 
$y_0 <0$.  Hence for all the corresponding points $D$ on the boundary of $\bar Q$, we have $E^2 \cdot D <0$. We show then that the linear function $E^2$ is negative on all of $ \bar Q$; if not, then for some point $D\in \bar Q$ we have $E^2 \cdot D =0$. 
  For such a point, the index condition on $E$ implies that $(D^2 \cdot E) E^3 \le (D\cdot E^2 )^2 =0$, and given that $E^3 >0$, this gives a contradiction when $D\in Q$. 
  More generally, if $E^2 \cdot D =0$ for a point of $\bar Q$, then we would deduce that $D\cdot E$ is numerically trivial and so $H(D)=0$.  Such points will be of the form $-D = D' = (0,y_1, 0)$ with $y_1 \le 0$, 
  where we have $E\cdot D'$ is explicitly the linear form 
  $-by_1 y - (1-b)(1-y_1)(1-x-y)$, whose term in $x$ has coefficient  $(1-b)(1-y_1)  <0$, contradicting the 
  previous deduction.
  
The conclusion is that $E^2$ is strictly negative on all of $\bar Q$.  In particular, if we choose an integral class $D \in Q$, then $h^0 (X, \O_X (mD)) >1$ for some explicit $m$ depending on $D$, but
 then for $\lambda \gg 0$, we will have that $(D - \lambda E)^2$ will be 
negative on $Q$, in particular on any ample class $L$, and such $D- \lambda E$ would not be movable.  Thus we can find a bound $\lambda _0$ for which $D- \lambda E$ is not 
movable for $\lambda > \lambda _0$.  Hence there is an upper bound on $r$ for which $mD - rE$ is movable, and so we find a finite number of classes, one of which 
will have a multiple which is mobile.
Using the argument used previously in the proofs of Propositions 3.3 and  4.1, a finite set of integral classes may be found, at least one of which is movable and big; thus 
  $X$ lies in a bounded family.
\end{proof}\end{prop}

Note that if we allow two rigid non-movable surfaces, then the cubic can be the Fermat cubic, as for instance is the case if we take a resolution of a quintic hypersurface with two simple nodes.
The Fermat cubic is the only example I know where an elliptic curve occurs  (defined by the ternary cubic  for a Calabi--Yau threefold with $\rho =3$) 
having only one component.  As we saw in 
Section 1, the Fermat cubic is unique in that any component of the positive index cone 
has boundary, part of which is curved contained in $F=0$, and the rest 
of which has two linear parts contained in $H=0$.

 We  now consider the remaining case, namely $k=-2$.  As observed in the Introduction, there will be a unique rigid non-movable surface  class $E$ on $X$.
 We assume without loss of generality that the boundary of the relevant component $P^\circ$ of the positive index cone corresponds to the affine branch of the cubic $C_1$ in the quadrant $x> 1/3$, $y>1/3$, and the appropriate line segment $C_3$ at infinity joining the inflexion points $B_1 = (0:1:0)$ and $B_2 = (1:0:0)$.  For a class $A$, we know when and where $G_A$ vanishes on $C_1$ and we explicitly calculate when and where it vanishes on $C_3$.  Unsurprisingly, given that it is a 
 limit of cases we know about 
 with $k \to -2 $ from above, $G_A$ is either non-zero on $C_1 \cup C_3$, or vanishes twice, except for the special case when $A= (1/3 , 1/3, 1)$, where $G_A (x,y,z) = -{1\over 3}z^2$.  For this point, 
 we would need to have $E$ a negative multiple of $A$, and since $A^3 = -1/3$ only a bounded number of multiples are allowed; choosing integral $D \in P^0 =Q$ we deduce as in 
 the proofs of Propositions 3.3 and 4.1 that  $X$ lies in a bounded family, 
 given that we can bound $\lambda$ with $H(D-\lambda E) \ge 0$.   Otherwise
 $E$ either represents a point in the affine plane $z=1$ other than $(1/3, 1/3, 1)$, or represents a point at infinity.
 In the former case it is clear that none of $C_3$ is visible from $E$, and so a contradiction may be obtained via Proposition 1.2.  So we assume that  $E$ 
 is on the line at infinity (in particular $H(E) =0$).  For $E$ on the line at infinity, we know that $E\cdot D^2$ defines a line pair intersecting at the affine point $(1/3:1/3:1)$, and that the  line $E^2 \cdot D$ is a further line through this point (these statements follow since they are the limit of facts we know when $k$ is near $-2$ or 
 from a direct calculation).  In particular, we may identify 
 these lines by calculating the relevant point at infinity. 
 Clearly $E$ will not correspond to a point in the closed cone generated by $(-1,0,0)$ and $(0,-1,0)$, since we know that for any ample $L$ we must have $E\cdot L^2 >0$, and it is not in the open cone generated by $(1,0,0)$ and $(0,1,0)$ by Lemma 0.3.  If say $E$ were a positive multiple of $(0,1,0)$ or $(1,0,0)$, then this 
 represents a rational 
 point of the real elliptic curve, which lies also on the line $c_2 =0$ by Lemma 0.3, contrary to the assumptions we have made,  Thus without loss of generality we may assume that $E$ is a positive multiple of 
 $(-1,\mu,0)$, with $\mu >0$.  We note in passing that $F(-1,\mu,0)  = -9\mu (\mu -1) \le 0$ if and only if $\mu \ge 1$.
 The following facts may then be checked:
 
 (1)  One of the lines of $E\cdot  D^2 =0$ defines a plane only intersecting $P$ at the origin in $\R^3$ and the other line (which may be explicitly calculated from $E$) defines a 
  plane which cuts $P^\circ$ in two; this determines $Q = \{ D\in P^\circ : E\cdot D^2 >0$\},  as the part whose boundary includes $B_2 = (1,0,0)$.
 
 (2)  The second line in (1) intersects the line at infinity in the point $(1: t(\mu ):0)$, where $$t(\mu ) = \mu -1 +( (\mu -1)^2 +\mu)^{1/2}.$$  Moreover the line 
 $E^2\cdot D =0$ does not intersect (and the linear form is negative on) the quadrant $x > 1/3$, $y > 1/3$ 
 if ${1\over 2} < \mu <2$; hence the form is negative on 
 the non-zero elements of  $P$; the line is the asymptote through $B_1$ if $\mu = {1\over 2}$ and the asymptote through $B_2$ if $\mu =2$.  For $0 <\mu \le 1/2$, we show that on the affine plane $z=1$ the form 
$E^2 \cdot D$ is strictly negative below the line $y-1/3 = s(\mu )(x -1/3)$, where $s(\mu ):= \mu (2-\mu)/(1-2 \mu)$.
Algebraic manipulations show that $s(\mu ) > t(\mu )$ in this range, and so we still have $E^2 \cdot D$ 
is strictly negative on the non-zero points in the closure of $Q$.

 (3) If we \it know \rm the class $E$, a positive multiple of say $(-1, \mu , 0)$ as above, we also know $Q$ and 
  can pick any integral class $D$ in $Q$; then $h^0 (X, \O_X (mD)) >1$ for some explicit $m$ depending on $D$.  If $0 < \mu < 2$, then 
$-2D^2\cdot E + \lambda D\cdot E^2 <0$ for all non-zero $D$ in the closure of $Q$ and all $\lambda >0$.  Hence 
 there exists $\lambda _0 >0$ such that for $\lambda \ge \lambda _0 >0$, we will have that $(D - \lambda E)^2$ is negative on all of $Q$ (and in particular on the \Kahler cone).
 Therefore there is a known upper-bound $\lambda _0$ on the $\lambda$ for which $D-\lambda E$ can be movable.  The case when $\mu >1$ is even easier; here $E^3 <0$  and we consider $f(t) = (D-tE)^3$, with 
 positive real roots $\lambda _0 > \lambda _1 >0$ say.  Suppose now that $L = D-\lambda E$ is mobile for 
 some $\lambda >\lambda _0$; then $L^3 >0$.  Any ample class $H$ in $P$ must have $z>0$, and then 
 $(L+ sH)^3 >0$ for all $s>0$; this implies that $L\in P$, which is clearly nonsense.  Thus in this case too, 
 we have found an upper-bound $\lambda_0 $ for the possible $\lambda >0$ for which 
 $D-\lambda E$ can be movable. 
 In both cases therefore, the data produces a finite set of classes 
 $mD -rE$, at least one of which will have a multiple which is mobile.  Using the argument used previously in the proofs of Propositions 3.3, 4.1 and 4.2, a finite set of integral classes may be found, at least one of which is movable and big; thus 
  $X$ lies in a bounded family.
  
 (4) We show that the class of $E$ may be assumed known up to finitely many possibilities.  
 If $c_2 \cdot E \le 0$, then by Propostion 2.2 of \cite{WilBd} there are only finitely many possible values for $c_2\cdot E$ and $E^3$.  
Suppose now $c_2 \cdot E >0$;  if $E^3 \ge 0$, then 
 by Proposition 2.2 of \cite{WilBd} again  there are still only finitely possibilities for $E^3$ and $c_2 \cdot E$.
 We may assume therefore that $E^3 <0$, and so some positive multiple corresponds to a point $(-1,\mu ,0)$ at 
 infinity with $\mu  >1$;
 from this it follows that $Q$ contains the open subcone of $P^\circ$ defined by $y <x$.  We can choose any integral $D$ in this smaller cone, and then for some explicit $m>0$, we have that $mD = \Delta + r E$ for some movable divisor $\Delta$  (some multiple of which is mobile) and some integer $r\ge 0$.  If $r=0$, 
 then we already have our desired big movable class, and as before we deduce boundedness and only finitely many possible values for $E^3$ and $c_2\cdot E$.   
  We may assume therefore that $r>0$, in which case the fact that $r c_2 \cdot E \le m c_2 \cdot D$ yields bounds on both $c_2 \cdot E$ and $r$.
 The upper-bound on $c_2 \cdot E$ will again yield only finitely many possible values for $c_2 \cdot E$ and  
 $E^3$ by Proposition 2.2 of \cite{WilBd}.  Under the assumption that the line $c_2=0$ does not meet the elliptic curve at an inflexion  point, we can deduce as in 
 the first part of the proof of Lemma 3.7 that there are 
 only finitely many possible classes for $E$.  Having  now reduced to the case when there are only finitely many possible classes for $E$,  we may run the argument from (3) to deduce the claimed result.
  
 \smallskip
 We have now dealt with the second case where the Hessian is singular.
 \begin{prop} Suppose that the Picard number $\rho (X) =3$ and the cubic is smooth with one real component  and singular Hessian, where 
  the line $c_2=0$ does not intersect the real elliptic curve at an inflexion point.
 If $X$ contains no more than one rigid non-movable surface, then there will exist such a surface and $X$ lies in a bounded family.
\end{prop}

\smallskip The Main Theorem from the Introduction now follows from combining
Proposition 2.4,  Corollary 3.8 and Propositions  4.1 and 4.3.

\bibliographystyle{ams}

\end{document}